\documentclass[letterpaper, 10pt, conference]{ieeeconf}
\IEEEoverridecommandlockouts

\newcommand{\Attempt}[1]{}  

\usepackage{wrapfig}

\usepackage{graphicx}
\usepackage{float}
\usepackage{subfigure}
\usepackage{derivative}
\usepackage[utf8]{inputenc}
 
 \newcommand{\RM}[1]{#1}

 \newcommand{\nchanges}[1]{{#1}}

\usepackage{amsmath,mathptmx}
\usepackage{amssymb}

\usepackage{dsfont}
\usepackage{graphicx}
\usepackage[colorlinks=true, allcolors=blue]{hyperref}

\newcommand{\eop}{{\hfill $\Box$}}
\newtheorem{thm}{{\bf Theorem}}
\newtheorem{lem}{{\bf Lemma}}

\newcommand{\ueta}{{\bar \eta}}
\newcommand{\leta}{{\underline \eta}}
\newcommand{\B}{B^{(2)}}
\newcommand{\ignore}[1]{}

\IEEEoverridecommandlockouts                              

\newcommand{\TR}[2]{#2} 

\usepackage{bbm}
\usepackage{enumerate}

\newcommand{\exactproof}[1]{}  

\begin{document}
\title{\bf Integrative Modeling and Analysis of the Interplay Between Epidemic and News Propagation  Processes
}

\author{Madhu Dhiman, Chen Peng, Veeraruna Kavitha, and Quanyan Zhu \thanks{M. Dhiman and V. Kavitha are with the Indian Institute of Technology, Bombay. Emails: \{214190003,vkavitha\}@iitb.ac.in; C. Peng and Q. Zhu are with the Department of Electrical and Computer Engineering, New York University. Emails: \{cp3509, qz494\}@nyu.edu.}}
\maketitle
\begin{abstract}
    The COVID-19 pandemic has witnessed the role of online social networks (OSNs) in the spread of infectious diseases. The rise in severity of the epidemic augments the need for proper guidelines, but also promotes the propagation of fake news-items. The popularity of a news-item can reshape the public health behaviors and affect the epidemic processes. There is a clear inter-dependency between the epidemic process and the spreading of news-items. This work creates an integrative framework to understand the interplay. We first develop a population-dependent `saturated branching process' to continually track the  propagation of   trending news-items on OSNs.
    A two-time scale dynamical system is obtained by integrating the news-propagation model with  SIRS epidemic model, to analyze  the holistic system. %
    %
    It is observed that a pattern of periodic infections emerges under a linear behavioral influence, which explains the waves of infection and reinfection that we have experienced in the pandemic. We use numerical experiments to corroborate the results and use Twitter and COVID-19 data-sets to recreate the historical infection curve using the integrative model. 
\end{abstract}

\section{Introduction}

Fake news becomes increasingly prevalent with the rise of online social networks (OSNs) \cite{benkler2018network,kapsikar2020controlling}. It is often generated for marketing or political purposes by social bots to manipulate public opinions and decisions. 
It's impact has gone beyond the influence on opinions and affected the behaviors of the people. COVID-19 is an example of this outcome;  for example during the pandemic, the fake news regarding vaccination has spread widely among social networks and caused an adverse effect on the prevention and mitigation of COVID-19 infection. Further, the deterioration of the pandemic  has in turn accelerated the creation and
propagation of fake news-items. 


Thus the spreading of fake news and the epidemic processes are interdependent.  
There is a need to understand the interplay between the two processes. There are several works discussing the relationship between information and epidemic propagation, including \cite{granell2014competing,liu2022herd,huang2022game,zhu2022preface}. They have focused on the impact of  decisions and awareness on the propagation of epidemic, assuming that individuals  take the same preventive measures once they receive the disease-related information. Very few works have explicitly modeled the dynamics of news propagation to scrutinize its inter-dependency with the spreading of infectious diseases. 



We propose an integrative model to capture and study the interplay between the propagation of trending news-items (fake as well as authentic) and the epidemic. 
We first develop population-dependent branching process with saturation, to model the propagation of news-items on OSN and their replacement by newer topics. 
This part is inspired by the latest work on total-population-size-dependent branching processes considered in \cite{agarwal2021new,agarwal2022new}.
The epidemic-related posts on OSNs can either get viral (i.e., a large number of the copies are shared) or become extinct. A news-item stops trending when its propagation saturates (i.e., it is no longer forwarded), and then is replaced by another news-item on the similar topic. We use  stochastic approximation based approach  to capture the dynamic life cycles of news propagation. 

Secondly, we use a Susceptible-Infectious-Recovered-Susceptible SIRS compartmental model to study the spreading of infectious diseases (e.g.,  \cite{lahrouz2011global}). The infection rate depends on the population behavior, which may be  influenced by the circulating news-items. We consolidate the news propagation model into the SIRS model to create an integrative framework that captures their mutual influences. We create a two-timescale dynamic process in which news-items spreads faster than the epidemic; this is justified by the fact that the online news-items evolve at a much faster rate than physical human contacts. 
Furthermore, the slower process, in this case, depends on the evolution of the faster process over a window of time, which calls for a   different type of two-time-scale process. 
The two-time-scale consolidated framework describes how news propagation can reinforce the spreading of epidemic (without interventions). 

We analyze the dynamical system properties of the coupled system under different behavioral influences by the news content. In particular, we observe that a pattern of periodic infections emerge under a linear behavioral influence, arising from the existence of limit cycles of the dynamical system. It explains the waves of infection and reinfection that we have experienced in the pandemic. We use Twitter and COVID-19 datasets to validate the proposed model. We successfully recreate the historical infection curve by fitting and identifying the timing of news influences across the 2-year period of the COVID-19 infection. 

\RM{
\noindent{\it Organization:} 
In section \ref{subsec_CP}, we introduce the OSN content propagation model. SIRS model  is considered in section \ref{subsec_epidemic}. In section \ref{sec:interplay}, we discuss the interplay between the two and analyze the dynamical system under different behavioral patterns. Section \ref{sec_numerical} presents numerical results.}



\section{News Propagation Model}\label{subsec_CP}
  Consider that there are $M$ epidemic-related tweets, posts, or news-items on the OSN. We are interested in tracking their influence, which can affect the ongoing epidemic.
When a post is 
 designed by a content provider, she shares it with an initial set of users called \textbf{seed users}. Each  recipient shares it with all or a subset of their friends, depending on the interest generated by the post. Let $\eta_m$ be the attractiveness factor of 
  $m$-th post, i.e.,  the probability that a typical recipient forwards the $m$-th   post to each of her friends, independently of others. A  friend of each recipient, upon reading the post, can again forward the post  (with probability $\eta_m$). The post propagation continues in this manner.

At any instance of time, among the copies forwarded/shared by then, some  recipients might not have read the post; we refer such copies as {\bf live copies}, while the rest are referred to as {\bf dead copies}.  Once a user reads the news-item, it is unlikely that she would be interested in the same news-item again. Thus, only the users with the live copies are responsible for sharing the post and continuing the propagation. {Let $\Psi_{m,k}$ be the number of live copies of the $m$-th news-item, immediately after  the $k$-th user forwards the post.} Further, let the total copies, including both live and dead ones, of the $m$-th news-item
be represented by $\Theta_{m,k}$. Then,  the propagation of the $m$-th  news-item can be captured as follows: for every $m \le M$,

\vspace{-3mm}
{\small\begin{equation}
 \Psi_{m, k+1} = \Psi_{m, k} - 1 + \xi_{m, k+1} \mbox{ and } \Theta_{m,k+1} = \Theta_{m,k} + \xi_{m, k+1} . 
 \label{Eqn_random_dynamics}
\end{equation}}
In the above,
  each  live copy is considered to be a \textit{parent}, which can generate multiple additional copies $\xi_{m,k+1}$ of the same topic, analogous  to a random number of offsprings in the branching process.

OSN users seldom read and forward the same news-item the second time. When a user forwards the news-item to her friends, some of them might have already received the post in the past. Thus the number of effective news forwards  (i.e., the offsprings) depends on the number of people who have already received a copy of the same news-item, i.e., the total copies. Hence, we require a population-size dependent branching process to model the news-propagation processes  accurately (as in   \cite{agarwal2022new,agarwal2021new}).
Furthermore, unlike the models considered in the literature, we require  the  branching processes whose offspring's depend on the total population. Such  branching processes have been considered recently in \cite{agarwal2022new,agarwal2021new}
and we consider a simpler modification of the same. 


We analyze the discrete-time dynamics corresponding to the  embedded chain,  the chain obtained by observing the   system immediately  after the transition epochs. To be precise, we study the ratios related to $\Psi_{m,k} \  \text{and} \ \Theta_{m,k}$ given in \eqref{Eqn_random_dynamics}. We refer to these epochs as \textbf{wake-up epochs}, as these are the instances at which a user  visits its timeline, reads  and forwards the news-items.
The influence of the top $M$ posts  (that affect the ongoing epidemic) is captured by the following ratios:
denote the ratio of the number of live copies  to the number of times the news-items is being forwarded  and that corresponding to the total  copies respectively by:
\begin{equation*}
    \psi_{m,k}  =\frac{\Psi_{m,k}}{k}, \ \mbox{ and }    \  
    \theta_{m,k}=\frac{\Theta_{m,k}}{k}, \mbox{ for any } k \ge 1.
\end{equation*}
 Recall that ${\xi}_{m, k+1}$ is  the effective number of friends to whom the $m$-th news-item is forwarded by  the $k$-th user. From \eqref{Eqn_random_dynamics}, it is clear that the above quantities can be rewritten in the following iterative manner with $\epsilon_k := {1}/{k}$:
\begin{equation}
    \begin{aligned}
    \label{Eqn_SA}
        \psi_{m, k+1} &= \psi_{m,k} 
+ \epsilon_{k+1}  \left ( 
\xi_{m, k+1} - 1 - \psi_{m,k}\right )
  \mbox{ }    \\
\theta_{m, k+1}& = \theta_{m,k} 
+ \epsilon_{k+1}    \left ( 
\xi_{m, k+1}   - \theta_{m,k}\right )
.
    \end{aligned}
\end{equation}
It takes the form of a stochastic approximation-based iterative scheme (\cite{benven}). It is important to note here that{ we ignore the time scale of the news-propagation system and that  the time scale does not play a role when analyzing \eqref{Eqn_random_dynamics} defined above.} 

\nchanges{
We now construct an appropriate stochastic iterative 
scheme in the following, that captures the news-propagation dynamics accurately using \eqref{Eqn_SA}-- it has to capture regular news propagation updates, replacement of a news-item by a newer one and continual tracking of  the news-items.  
}

There are two kinds of updates related to $\Psi_{m,k}$. The first  kind of update arises from the regular news  propagation, i.e., driven by user forwards.  The second kind of update occurs  when a   news-item stops spreading and a fresh news-item (of the similar topic) emerges.
It happens when the effective-forwards related to an old news-item diminish as the total-copies reach a saturation level;  as the number of forwards  $\{\xi_{m,k}\}$ diminish, the live-copies $\{\Psi_{m,k}\}$ reduce (see \eqref{Eqn_random_dynamics}), while the total copies $\{\Theta_{m,k}\}$ saturate; thus saturation is reached at $k$-th epoch  if $\Psi_{m,k} \le \delta_\psi$ and $\Theta_{m,k} \ge  \delta_\theta$, for some appropriate $\delta_\theta, \delta_\psi$. Hence the saturation point is captured by regime $\mathcal{E}^2_{m,k}$ given below. Here, $ \mathcal{E}^1_{m,k}$ represents the first regular news-item-update regime.

\vspace{-4mm}
{\small\begin{eqnarray*}
    \mathcal{E}^j_{m,k} &=&\{(\psi_{m,k},\theta_{m,k})\in\mathcal{R}^j\},\  \text {for}  \  j = 1,2, \text{ where} 
\\
\mathcal{R}^2 &=&\{(\psi, \theta):0\leq\psi\leq \delta_\psi \text{ and }  \theta\geq \delta_\theta\}, \  \ 
  \ \mathcal{R}^1 =  (\mathcal{R}^2)^c 
 .
\end{eqnarray*}}
  We are interested in the tracking   performance, and thus $\epsilon_k$ has to be replaced with a constant $\epsilon > 0$ in \eqref{Eqn_SA} (e.g., \cite{benven}).   Towards saturated news-item replacement,   we introduce extra or fictitious iterates that  replace the news-item with a new one\footnote{We assume that the social network supports a number of simultaneous posts, and there is always a newer post that can replace a just saturated one.};  once an old news-item saturates (say $m$), i.e.,    
${1}_{{\mathcal E}^2_{m,k}} = 1$,    the additional iterates (using constant  $-C$ in \eqref{Eqn_update_stocahstic} given below) keep reducing the corresponding  total population till it reaches close to zero (i.e., below $\delta_\theta$).  
Hence the overall update  equation is given by (for each $m$):

\vspace{-4mm}
{\small \begin{eqnarray}    \label{Eqn_update_stocahstic}
        \psi_{m, k+1} = \psi_{m,k} 
+  \epsilon \left( \mathcal{J}_{m,k} \left ( 
\xi_{m, k+1} - 1 - \psi_{m,k}\right )- \mathcal{J}^{'}_{m,k}( \psi_{m,k})\right) 
, \hspace{-26mm} & \nonumber \\
\theta_{m, k+1}  = \theta_{m,k} 
+ \epsilon\left(  \mathcal{J}_{m,k} \left ( 
\xi_{m, k+1}   - \theta_{m,k}\right ) - \mathcal{J}^{'}_{m,k}(C\theta_{m,k})\right), \hspace{-20mm} &\nonumber \\
  & \hspace{-46mm} \mbox{\normalsize where } \mathcal{J}_{m,k} = {1}_{{\mathcal E}^1_{m,k}} \mbox{\normalsize and }  \mathcal{J}^{'}_{m,k} = {1}_{{\mathcal E}^2_{m,k}}.\hspace{7mm}
\end{eqnarray}}
\ignore{
{\color{red}
From \eqref{Eqn_update_stocahstic}, there can be correlations between the propagation processes of two concurrent news-items either when they   compete  or occupy the timeline of the same user. In this work, we consider two decoupled news-items, which generate independent interests among the users. We also ignore the effect of  timelines (and users' lethargy to read all the content displayed on their timelines) on news propagation. It is clear  that one can analyze the news propagation corresponding to each of the $M$-trending topics separately. 
}}
\vspace{-9mm}
\subsection*{ODE approximation}
The dynamics of news-items related to  various channels/trending-topics is independent of each other. Hence, it is sufficient to
analyze the propagation of each post  independently of the others. Fix any $m$. 
Let ${\cal F}_k :=\sigma \left ( (\psi_{m,s}, \theta_{m,s}) : s \le k  \right )$ be the sigma algebra generated by the $m$-th post dynamics  till epoch $k$. Observe that the corresponding conditional expectation is given by:

\vspace{-3.5mm}
{\small
\begin{eqnarray}
E_k \left[  \mathcal{J}_{m,k} \left ( 
\xi_{m, k+1} - 1 - \psi_{m,k}\right ) + \mathcal{J}^{'}_{m,k}(- \psi_{m,k})\right]  \\ =  \mathcal{J}_{m,k}\bigg( {\cal M}(\theta_{m,k} )-1-{\psi}_{m,k}\bigg )+  \mathcal{J}^{'}_{m,k}(-{\psi}_{m,k} ), 
\label{conditional_expectation}
\end{eqnarray}}
where $E_k[\xi_{m,k+1}] := \mathcal {M} (\theta_{m,k}).  $ 
We immediately have the following ODE approximation result applying the   theory of stochastic approximation under the   following assumptions:

 \begin{enumerate}[{\bf A}.1]
     \item \it For any $\lambda \in [0.5,1],$ there exists a constant $\varrho$ such that
     $|{
     \cal M} (\theta)   -  {
     \cal M} (\theta')| \le \varrho |\theta - \theta'|^\lambda.$
     
     \item 
     $P (\xi_{k+1} \le \rho | \theta_k) = 1$ for any $\theta_k$ and for  some $\rho < \infty.$
     
 \end{enumerate}
\begin{thm}
\label{newsode}
Assume {\bf A}.1-{\bf A}.2.  Let
$\phi^{(\epsilon)}_{m,k} := (\psi_{m,k}, \theta_{m,k} ) $
represent   the news-propagation trajectory   \eqref{Eqn_update_stocahstic}    at $\epsilon$ and let 
$\phi(\tau) := (\psi(\tau), \theta(\tau))$  be the solution of following ODEs:

\vspace{-4mm}
 {\small\begin{eqnarray}
        {\stackrel{\bullet} {\psi}}_m(\tau)
        \hspace{-1mm}
        & \hspace{-1mm}= \hspace{-1mm}&  \hspace{-1mm}\mathcal{J}_{m} (\phi)\bigg( {
     \cal M} (\theta_m)-1-{\psi}_m(\tau)\bigg )+ \mathcal{J}^{'}_{m} (\phi) (-{\psi}_m(\tau)),\nonumber \\
     {\stackrel{\bullet}{\theta}}_{m}(\tau)    \hspace{-1mm}
        & \hspace{-1mm}= \hspace{-1mm}&  \hspace{-1mm} \mathcal{J}_{m} (\phi) \bigg(  {
     \cal M} (\theta_m)-{\theta}_{m}(\tau)\bigg )+ \mathcal{J}^{'}_{m} (\phi) (-C{\theta}_m(\tau)), \nonumber \\
      \mathcal{J}'_{m} (\phi)     \hspace{-1mm}
        & \hspace{-1mm}= \hspace{-1mm}&  \hspace{-1mm} 1_{\{ \theta > \delta_\theta , \ \psi < \delta_\psi \} }, \mbox{ and, }\ \mathcal{J}_{m} (\phi) = 1- \mathcal{J}'_{m} (\phi).     \label{Eqn_update_ODE}
\end{eqnarray}}
For any $m$, as $\epsilon\rightarrow 0$, trajectory   \eqref{Eqn_update_stocahstic} converges to the solution of  above ODE in the following sense: for any  finite time $T$,
    \begin{eqnarray}
    P \left ( \sup_{k \le T/\epsilon}  \left | \phi^{(\epsilon)}_{m,k} -  \phi  (\epsilon k)   \right | > \delta \right ) \to 0 \mbox{ as } \epsilon \to 0,
    \end{eqnarray}
when the initial conditions are equal, $\phi(0) = \phi^{(\epsilon)}_0.$ 
\end{thm}
{\bf Proof}    \TR{is omitted due to  page limit (available in \cite{TR})}{is in Appendix}. \eop

 To capture the saturation effect because of re-forwards   discussed earlier, we consider the following linear model:
 \vspace{-4mm}
\begin{eqnarray}
\mathcal {M} (\theta) 
= \eta - a \eta \theta
 = \eta ( 1 - a \theta),
\label{Eqn_Mean_offsprings}
\end{eqnarray}
where   $( 1 - a \theta)$ represents the reduction in the expected number of effective forwards. \textit{The   conversion factor $a$
is clearly specific to an OSN, and it is the same  for all the posts}.

\subsection*{ODE Solutions  and the influence of a viral post}
From  \eqref{Eqn_update_ODE}, it is not difficult to show that $\psi_m(t) \rightarrow 0$; i.e., the post vanishes immediately, when ${\cal M}(0) = \eta <1 $. The post explodes and gets viral (i.e., the number of copies grow significantly), only when ${\cal M}(0) >1 $.  We derive the  solution of the ODE  \eqref{Eqn_update_ODE} in the following, to capture its possible influence on the epidemic. 
The RHS of the ODE \eqref{Eqn_update_ODE} is piece-wise continuous (in two regimes), and hence the solution is in the extended sense (i.e., it satisfies the ODE for almost all $\tau$). We first derive the  solutions in individual regimes,  connect them together by appropriate initial conditions, and then derive the asymptotic behavior of the solution. 

Suppose that regular update regime (while ${\cal J}(\phi) = 1$) is concluded at $\tau_1$; this happens when $\psi$ goes below $\delta_\psi$. Thus the initial conditions for the 
  next (replacement) regime are, {\small$\psi(0) = \delta_\psi$} and {\small$\theta(0)= \theta(\tau_1)$}. 
Now, 
    {\small $\dot{\psi}(\tau) = -{\psi}(\tau)$} and {\small$\dot{\theta}(\tau) = -C{\theta}(\tau)$.} 
Thus, the solution is,   {\small $ \psi(\tau) = \delta_\psi e^{-\tau}$} and {\small $ \theta(\tau) = \theta(0) e^{-C\tau}.$}

The next  regime   starts with {\small$\theta(0) = \delta_\theta$} and {\small$\psi(0) = \psi (\tau_2)$}, as the replacement regime is concluded  when $\theta$ drops below $\delta_\theta$ (and say at $\tau_2$).
From \eqref{Eqn_update_ODE} and \eqref{Eqn_Mean_offsprings}, 
{\small$\dot{\psi}(\tau) = \eta - a \eta  \theta(\tau) - 1 -{\psi}(\tau)$} and 
{\small$\dot{\theta}(\tau) = \eta - a \eta  \theta(\tau) -C{\theta}(\tau)$}. By solving these equations,

\vspace{-4mm}
{\small\begin{eqnarray*}
\psi(\tau) &=&  -\kappa_1 e^{-(a \eta +1)\tau} + \kappa_2 e^{-\tau} + \kappa_3, \  \kappa_1 = \eta -(a \eta +1)\delta_\theta, \\
\kappa_2 &=& \psi({0}) + a \eta ^2 - (a \eta +1 )\delta_\theta + 1 ,  \mbox{ and }
 \kappa_3 = \eta  - 1 - a \eta ^2,  \\
 \theta (\tau) &=& \left(  \delta_\theta - \frac{\eta}{a \eta +1 }   \right) e^{-(a \eta +1)\tau} + \frac{\eta}{a \eta +1 }.  
\end{eqnarray*}}
This solution converges to a limit cycle alternating between  two regimes (see Fig.\ref{fig:limit_cycle_News}). By solving appropriate fixed-point equations, with  initial conditions {\small$(\delta_\psi, \theta(0))$} for regime 1 and {\small$(\psi(0), \delta_\theta)$} for regime 2, we can obtain the limit cycle at the limit (see details  in the  Appendix). 

We assume  that \emph{the epidemic is influenced by news-propagation over  a window of time, and we capture this influence via the maximum value attained by the total population fraction at limit cycle}, which is given by 
equation \eqref{Eqn_istar_rstar} 
in the Appendix: \vspace{-1mm}
\begin{eqnarray}
\hspace{5mm}
\label{limit_of_frac}
  {\theta}_\infty ^{*} (\eta) := \lim_{k \to \infty} \max_{\tau' \in L_{C_k} } \theta (\tau') \approx \frac{\eta } { a \eta  + 1},
\end{eqnarray}where $L_{C_k}$ is the $k$-th  cycle.
 
\ignore{

we have two different scenarios with initial conditions $\psi(0) = \psi_{0}$ and $\theta(0) = \theta_{0}$ :

\begin{enumerate} 
    \item {\bf Phase (a):} 
    In the first region, i.e., $\mathcal{R}^1$, $\dot{\psi}(\tau) = -{\psi}(\tau)$ and $\dot{\theta}(\tau) = -C{\theta}(\tau)$. \\ 
    On solving above equations, we obtain $ \psi(\tau) = \psi_{0} e^{-\tau}$ and $ \theta(\tau) = \theta_{0} e^{-C\tau}.$ 
    \item {\bf Phase (b):} Consider first region, i.e., $\mathcal{R}^2$.
    As $$
{\cal M} (\theta) = \eta - a \eta  \theta$$ and in this region, \\
$\dot{\psi}(\tau) = \eta - a \eta  \theta(\tau) - 1 -{\psi}(\tau)$ and 

$\dot{\theta}(\tau) = \eta - a \eta  \theta(\tau) -C{\theta}(\tau)$.

On solving above equations, \\
$ \psi(\tau) = -\kappa_1 e^{-(a \eta +1)\tau} + \kappa_2 e^{-\tau} + \kappa_3$ and \\
$ \theta (\tau) = \left(  \delta_\theta - \frac{\eta}{a \eta +1 }   \right) e^{-(a \eta +1)\tau} + \frac{\eta}{a \eta +1 }  $, \\
where because of restrictions on $\mathcal{R}^2$ ,\\ $ \kappa_1 = \eta -(a \eta +1)\delta_\theta ,\\ 
\kappa_2 = \psi_{0} + a (\eta ^2 - (a \eta +1 )\delta_\theta ,  $ and \\
$ \kappa_3 = \eta  - 1 - a (\eta ^2).$

Therefore, we arrive at
\begin{eqnarray}
\label{limit_of_frac}
\lim_{\tau \to\infty}\psi(\tau) = 0 \mbox{ and } \lim_{\tau\to\infty}\theta(\tau) = \frac{\eta } { a \eta  + 1}.
\end{eqnarray}
\end{enumerate}
}

\section{Epidemics Model: Disease propagation}

 \label{subsec_epidemic}
We consider an epidemic population with three types of sub-populations,  the  infected (I), the  recovered (R), and the susceptible (S). 
People  who can be infected when in contact with any infected individual, are the susceptible sub-population, while people who recover from the infection and are immune to the disease are members of the recovered sub-population (remain recovered till they lose immunity).

 We consider a  SIRS compartmental  model that captures the dynamics of the relevant sub-populations. Let $i(t)$, $s(t)$, and $ r(t) $ be the normalized size of the infected,  the susceptible and the recovered population, respectively, at time $t$. Note that $i(t)+s(t)+ r(t)=1$, for all $t\geq0$.

The infectious disease spreads at rate $\bar \beta $, while the infected get recovered at rate $\alpha$. {In addition, we consider that a fraction $p_r$ of the recovered  is immunized, while the recovered (and the immunized) sub-population   can lose immunity at a rate  $l_i$.}  Thus the  dynamics are described as follows:
\begin{eqnarray}
\label{epi_ode}
    \frac{di}{dt} \ =  \ i \big ( \bar \beta(1-i - r)-  \alpha \big ),   
    \  \  \   \ 
\frac{d r}{d t} \ = \  \left  (   i \alpha p_r - r l_i   \right ).  
\end{eqnarray}

 In the following, we obtain the 
 asymptotic behavior of the above dynamics (proof is in Appendix) for any initial condition  (i.c.) $(i_0, r_0) \in \B := \{ (i,r) \in [0,1]^2: i+r \le 1 \}$. 
\begin{thm}
\label{thm_sir_attractor}
 The sub-population sizes $(i(t), r(t))$ given by \eqref{epi_ode}  converge to a unique limit for any i.c., in $\B$:

\begin{enumerate}[(i)]
    
    \item When $\bar \beta <   \alpha$, then the disease is eradicated eventually,  i.e., $(i(t), r(t)) \to (0, 0)$ as $t \to \infty$;
    
    \item When $\bar \beta >   \alpha$, the disease is not eradicated, i.e., $(i(t), r(t)) \to (i^*, r^*)$, where
    
    \vspace{-4mm}
    {\small 
    \begin{eqnarray}
    \hspace{-5mm}
i^* = \frac{\bar \beta -   \alpha}{\bar \beta c^r_i}  \mbox{ and } r^* = \frac{ \alpha p_r}{l_i} i^*
,\mbox{ where, } c^r_i := 1+ \frac{\alpha p_r }{ l_i}. \mbox{ \eop } \hspace{1mm} 
\label{Eqn_Only_Epidemics}
\end{eqnarray}}
\end{enumerate}
\end{thm}
Thus the sub-population sizes converge to a unique vector (disease-free if the limit is $(0,0)$)   when the dynamics are driven only by the epidemic.
However, this is not always the case under the influence of circulating news-items.  We will see the existence of  limit cycles, multiple limits,  etc.

\vspace{-2mm}
\section{Integrative Two-Timescale  System}\label{sec:interplay}
 The infection rate $\bar \beta$ in \eqref{epi_ode} represents the number of susceptible individuals who become infected by one infected individual per unit time (e.g., \cite{singh2021evolutionary}). 
 This rate is considered as a constant in classical SIRS models.
In reality,   this rate  is time-varying and is strongly dependent on human activities.  For example, it depends upon the  frequency of contacts and   the precautions taken by the individuals.  Furthermore,    $\beta(t)$ (rate at time $t$)  heavily relies on  the   information available  at time $t$. For example,   the viral news-items that  spread misinformation about the precautionary measures (e.g., social distancing, mask wearing, or  the ones that downplay the severity  of disease), can result in a higher $\beta(t)$ (as they promote reckless behaviors).  Likewise, the health policies from the government agencies  (e.g., recommendations to  stay   home) can make people 
  more cautious and  lower $\beta(t)$. 


On the other hand,  an increase in the infection level can breed panic in the population. As a result, people will seek more proactively for relevant news and information. Hence it is clear that the attractiveness factor   $\eta$ of \eqref{Eqn_Mean_offsprings} can depend on $i(t)$, the size of the infected subpopulation or the infection level at time $t$. The main aim of this paper is to study these interactions and analyze the resultant asymptotic outcome.

\ignore{
{\color{red}
In practice, the lifespan of news-items (hours/days/weeks) is typically shorter than that of the epidemic (months/years).
Thus to study the interplay one would require some type of two time scale process (e.g., \cite{borkar2009stochastic}). However   %
%
the propagation of epidemic is more aptly influenced by the behaviour of ensemble of people. The more the people that do not follow the guidelines (e.g., masks, social distancing  etc) the more the spread. And various individuals absorb these news-items at various instances of time. 
Thus the epidemic (evolving in slower time scale) is not influenced by news-items (evolving at faster time scale) at a given time instance, as considered in usual two-scale algorithms, but rather is influenced by happenings in some window of time. \emph{This calls for a different type of two-time-scale analysis.}
}}

The life cycle of news-items is often significantly shorter than that of an infectious disease.  The epidemic is more aptly influenced by the behaviour of  people who seek the news and information on OSN over a period of time.
To model such interactions, we consider a different type of two-time-scale system by connecting the two  dynamical processes  using ODEs; the (slow) epidemic-ODE is influenced by limit cycles \eqref{limit_of_frac} of (fast) news-propagation ODE as explained below. 
  Let $\tau$ be the time index for the fast timescale of news propagation and $t$ be the one for epidemic.

 We assume that the influence on epidemic at any time instance depends upon the outcome of the just concluded trending topics (e.g., the viral news-items).  As the news-items propagate at a faster time-scale, one can 
 capture the above influence using the maximum of the total population fraction at limit cycle, ${\theta}_\infty^{*} (\eta) $ given in \eqref{limit_of_frac}  of the latest trending topics. Note that we allow this fast time-limit to depend on the state of the epidemic at time $t$, via $\eta$. 
 
Varied interest toward news-items can  
 %
 %
%
%
be captured by modelling that $\eta_m$,  and hence, the  expected effective forwards depend on  $i$ 
(see \eqref{Eqn_Mean_offsprings}),
$
 {\cal M}_m (\theta; i) = \eta_m(i) - a \eta_m (i) \theta. $

 \ignore{\color{red}
People can respond to the consumed content differently depending upon the specifics of the content and the level of the infection prevailing in the area. For example, if content is predominantly  fake and is spreading wrong information about the disease, people can become irresponsible and the disease spread can get further accentuated.  On the other hand, if the  consumed content is dominated by authentic information published by responsible authorities, people can become more cautious and one can probably diminish the rate of disease spread. }

The response of the population clearly depends on the content of the news. In addition, it  depends on the level of  infection prevailing in the area. 
%
This aspect is captured via sensitivity parameters $\{w_m(i)\}$; in all, the infection rate is added with an amount, $w_m(i) {\theta}_\infty ^{*} (\eta(i)) $, due to the propagation of $m$-th news-item.
%
%
%
For a post that spreads mis-information,  $w_m(i)$  is  positive, it   is negative  for   authentic content. 
The overall rate of infection at time   $t$ is (see \eqref{limit_of_frac}):

\vspace{-4mm}
{\small\begin{eqnarray}
\label{Eqn_beta_t}
\beta(t)= \beta(i(t))&=& \bar \beta + \sum_{k=1}^M w_m(i(t)) {\theta}_\infty ^{*} (\eta(i(t)) ),\nonumber \\
&= &\bar \beta + \sum_{m=1}^M 
w_m (i(t)) \frac{\eta_m (i(t))  } { a \eta_m (i(t))  + 1} .
\end{eqnarray}} 
As observed above, \emph{$\beta(t)$ is time dependent only via  $i(t)$ and henceforth  we refer it as $\beta(i)$.}
{We assume    that the news-item attraction factors are of  two levels, i.e.,      $\eta_m \in \{ \underline {\eta}, {\bar \eta}\}$}   for any $m$. Either  a   particular news-item 
generates a lot of interest  and becomes viral (i.e., when $\eta_m =  \ueta$), or it gets extinct  without making an impact  (i.e., when   $\eta_m = \underline {\eta}$). 
Thus the influence of a post  is nonzero,  is captured by \eqref{limit_of_frac}, only when $\eta_m =  \ueta$. So, 

\vspace{-9mm}
{\small\begin{eqnarray} \hspace{10mm}
 {\bar w} (t) := \sum_{m=1, \eta_m = \ueta }^M 
w_m (i(t))  = {\bar w}(i(t)), \label{Eqn_barw_i}
\end{eqnarray}}is the (infection-level-dependent) integrated influence of all the viral posts at a given instance of time. {We assume that 
 this factor is deterministic\footnote{The number of trending topics, $M$,  is large, and one can justify this assumption using the law of large numbers.}, and it depends only on  $i$.} It is possible that fake and authentic news-items can co-exist and co-propagate at the same instance of time, with the former adversely influencing the disease propagation,  and $
{\bar w} (i)$ in \eqref{Eqn_barw_i} is the overall influence.

By  
incorporating  the above-mentioned factors into  \eqref{epi_ode}, we obtain a consolidated ODE  that captures the interplay between the two processes:

\vspace{-4mm}
{\small\begin{eqnarray}
\frac{d i}{d t} &=& i \left  (  \beta (i) (1-i-r) -  \alpha \right )   \label{Eqn_combined_ODE}\\
\frac{d r}{d t} &=&  \left  (   i \alpha p_r - r l_i\right ) \mbox{, with, }  
\beta(i) \ = \ \bar \beta +   \bar{w}  (i) \frac{\eta (i)  } { a \eta (i)  + 1}. \nonumber
\end{eqnarray}}
 This  consolidated ODE is instrumental in analyzing   a variety of scenarios.  We present a number of important scenarios in the following.
 We begin with a  case in which the population's interest in news-items increases with infection level.  Prior to this, we provide a few definitions. 

 \noindent{\it Asymptotic behavior:}
 We would analyze the consolidated ODE \eqref{Eqn_combined_ODE} to understand the time-limiting behavior using the results of two-dimensional ODEs. We observe different types of asymptotic behavior,  when the dynamics start in 
 ${B}^{(2)}$, i.e., for any initial condition (i.c.), 
 $(i_0, r_0) \in {B}^{(2)}$:

{\bf Local attractor (LA):} We refer a point $(i^*,r^*)$  as a local attractor, if  there exists a neighbourhood $\mathcal{N}$ such that $(i(t),r(t))\to (i^*,r^*)$ for all i.c. $(i_0,r_0) \in \mathcal{N}$.

{\bf Global Point Attractors (GPA)} In this regime, as $t \to \infty$, the solution of ODE $(i(t), r(t) )$ converges to an equilibrium  point in $B^{(2)}$, depending upon the i.c.; at maximum there are two equilibrium points, one of them is $(0,0)$.

{\bf Closed Orbits or Point limits (CoP)} Here, $(i(t), r(t) )$   converges to  one of the two point limits (one of them is disease free and a saddle point;   the other  is an LA) or to a closed orbit (limit cycle), depending upon the  i.c. 

{\bf Predominantly Closed orbits (PCo):}  
For some  set of  i.c.s, $(i(t), r(t))$ converges to disease-free state $(0,0)$, which is a saddle point. For the rest, the dynamics  either converge to  a closed orbit around an unstable equilibrium point, or to\footnote{We strongly believe it does not converge to the unstable point (also well understood in literature, e.g., \cite{verhulst2006nonlinear}), but require certain technical conditions to complete this proof and are on the way to proving the same. } the unstable equilibrium point itself.

 \subsection{Increasing Interest in News (I3N)}
 
     Interest in reading and forwarding relevant news or information increases with an increase in infection level, $i$. In this case,  $\eta(i)$ increases with $i$,  and hence we let   attractiveness factor take a linear form, i.e.,  $ \eta  (i) = \bar \eta (p i  + q ) $, with $p, q \ge 0$. We let the response to the news independent of $i$; i.e.,   ${\bar w} (i) \equiv {\bar w}$.
The  consolidated model  \eqref{Eqn_combined_ODE} for this case is:

\vspace{-4mm}
{\small 
\begin{eqnarray}
\label{Eqn_linear_influence_ode}
 \frac{di}{dt} = i \left ( \left (\bar \beta + {\bar w} \frac{\bar \eta (p i + q) }{ a \bar \eta (pi+q) + 1} \right ) (1-i-r) -   \alpha \right ).
 \end{eqnarray}}
 Let $\beta_0 := \beta (0) = \bar \beta + \bar w \bar \eta q / (a \bar \eta q + 1)$ represent the infection rate at $i=0$. 
 The asymptotic analysis of this scenario depends upon $\beta_0$ as given below  (see the proof in the Appendix): 
 \begin{thm}
\label{thm_Increased_interest} Consider the I3N case given in \eqref{Eqn_linear_influence_ode}.
 \\
  (i) When $\beta_0 <   \alpha$, the disease-free state (0,0) is an LA. If in addition, 
   either  $\bar w< 0$ or    $\bar\beta a + {\bar w} - \alpha a < 0$, the dynamics settle to disease-free  GPA regime  with $(0,0)$ as the only limit. 
   \\
(ii) 
     When $\beta_0 >   \alpha$,  the dynamics settles to the CoP regime with $(0,0)$ and $(i^*, r^*)$ as  the  two possible limits, where
  
    \vspace{-4mm}
    {\small 
    \begin{eqnarray*}
i^* &=& \frac{ c_b + \sqrt{ c_b^2 - 4 c_a c_c}}{-2 c_a}, \  r^* = \frac{\alpha p_r}{l_i} i^*\mbox{ \normalsize  with, }  \\
&& \hspace{-13mm}
c_b =  (\bar\beta a +{\bar w} -{\alpha} a) \ueta p - c^r_i(\bar\beta a+{\bar w})\ueta q - c^r_i \bar\beta  , \\
&& \hspace{-13mm}
c_a = - c^r_i (\bar\beta a + {\bar w} ) \ueta p, \mbox{\normalsize and, } 
c_c = (\bar\beta a +{\bar w} -{\alpha}a) \ueta q + \bar\beta -{\alpha} . \mbox{  \eop}
\end{eqnarray*}}
\end{thm} 
\noindent {\bf Remarks:} When authentic posts get viral ($\bar w < 0$), the disease can be curbed more easily: consider that $\beta_0 < \alpha$, but that  $\alpha < \bar \beta$;  then by Theorem \ref{thm_sir_attractor}, the disease is not cured under SIRS dynamics; however, it is curbed in presence of authentic news as confirmed by part (i). 

Similarly, with fake posts spreading misleading information, the diseases can prolong further (e.g.,  if  $\beta_0 > \alpha > \bar \beta$). 

The limit infection levels  are different with and without considering the influence of content.  
 
 In part (ii), the dynamics settle to the CoP regime, and here $(i^*,r^*)$ is an LA. Technically one can also converge to  a limit-cycle or $(0,0)$ (which is a saddle point);  however, we notice through several numerical examples  \TR{(see \cite{TR})}{(see Fig. \ref{fig:limit})} that the dynamics always converge to $(i^*,r^*)$. We do not observe periodic behaviors  for the I3N scenario.  
 \TR{}{
 \begin{figure}[htbp]
    \centering
  \vspace{-2mm}   \includegraphics[scale=0.2]{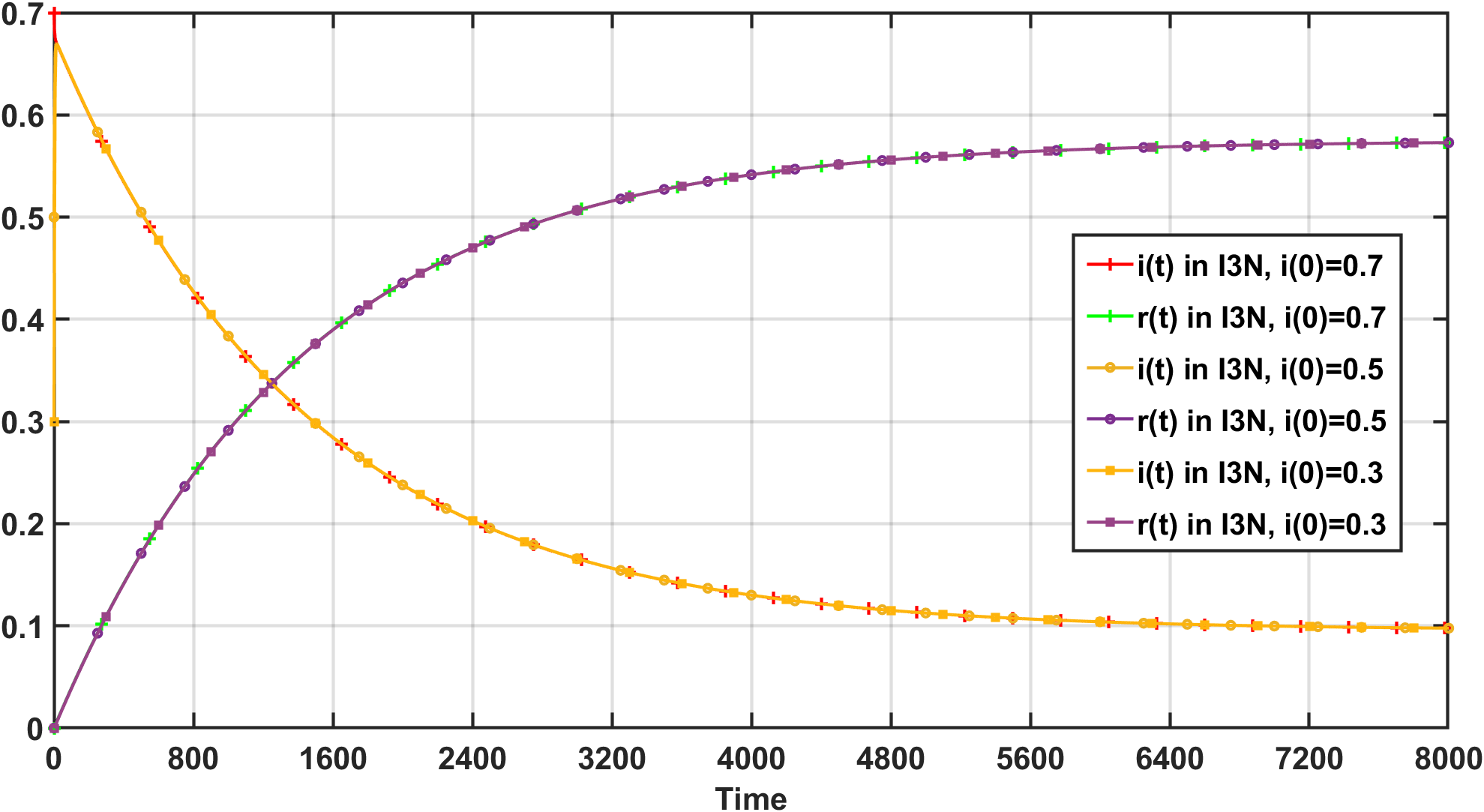}
  \vspace{-3mm}  \caption{Trajectory settles to $(i^*,r^*)$}
    \label{fig:limit}\vspace{-1mm}
\end{figure}}

 \ignore{
 \begin{thm}
\label{sir_linear_beta}
i)
When $\beta < \alpha(1-p_r)$ and when
$..$ then the disease is eradicated, i.e.,  $(i(t), r(t)) \to (0, 0)$.

{\color{blue}
$ \hat{C}>>  \alpha k(1-p_r)$ then  the disease is not eradicated but settles to unique value, i.e., $i(t) \to i^*$. Also 
\begin{eqnarray}
i^* = \frac{(\hat{C}-\beta k) +  \sqrt{(\beta k - \hat{C})^2 + 4\hat{C}k(\beta - \alpha) }}{2 \hat{C}k} \\ \text {and}  \ 
r^* = \frac{i^* {\tilde\alpha} p_r}{l_r}, 
\end{eqnarray}

where $ k = \left( 1+ \frac{\alpha p_r}{l_i} \right)$. Otherwise $(i(t) , r(t)) \to (0, 0)$}
\\
ii) When $\beta > \alpha(1-p_r)$, the disease is not eradicated and we have $(i(t), r(i) ) \to (i^*, r^*)$ where:

{\bf Proof:} See Appendix C.
\end{thm}}
 
 \subsection{
Increasing Behavioral Influence by News (IBIN)
}

When an individual reads news about the epidemic, the response can depend on the infection level. For example,  posts that promote mask wearing can positively influence individuals by encouraging them to  wear masks  when the infection level $i(t)$ is high. The same  individual may not follow the recommendation when the infection level is small. Similarly, the response to fake news can be different. Thus the behaviors of the individuals can depend on the infection level. 
Here, we consider a scenario, in which, the population's interest in news-items is constant, but their behaviors are influenced by the news. 
Such scenarios can be captured by
letting $\eta(i) \equiv 1$ (and $a=0$) 
and with
$
  \bar{w} (i) =   
 u  i $ (linear influence).    The constant $u$ is negative when the news-items are predominantly authentic; it is positive otherwise. 
The consolidated ODE \eqref{Eqn_combined_ODE} in this scenario is given by:

\vspace{-3mm}
{\small\begin{eqnarray}
\frac{di}{dt} = i \big ( \left ( \bar\beta +  u i \right ) (1-i-r) -  \alpha \big ), \mbox{ and} 
   \ 
\frac{d r}{d t}  =   \left  (   i \alpha p_r - r l_i\right ).
\label{Eqn_linear_w}
\end{eqnarray}}
We have the following result with proof   in the Appendix.
 \begin{thm}{\bf [Limit-Cycle and Reinfections]}
\label{sir_linear_beta}
Consider the consolidated ODE (\ref{Eqn_linear_w}) under the IBIN scenario. \\(i) When $\bar\beta <  \alpha$, $(0,0)$ is an LA. Further if $u < \alpha  c_i^r$,
    we have the disease-free GPA regime with 
    $(0, 0)$ as the only  limit. \\
(ii) When $\bar\beta >   \alpha$,   the disease need not be eradicated and, 
\begin{enumerate}[(a)]
    \item  If $u(1-i^*-c_i^r i^* )\le \bar\beta +   l_i/{i^*}$,   we   have  the CoP regime with possible limits from $\{ (0,0), (i^*, r^*) \}$ where:

\vspace{-4mm}
{\small 
\begin{eqnarray}
\hspace{-3mm}
i^* = \frac{(u  -\bar \beta c^r_i) +  \sqrt{(\bar\beta c^r_i - u )^2 + 4u  c^r_i(\bar\beta - \alpha) }}{2 u  c^r_i} \nonumber   \text { and}  \ \ 
r^*  = \frac{i^* {\alpha} p_r}{l_i}, 
\end{eqnarray}}
\item  If $u(1-i^*-c_i^r i^* ) >  \bar\beta +   l_i/{i^*}$ we have the PCo regime, where  $(i(t), r(t))$ eventually follows a periodic path or settles to $(0,0)$  or $(i^*,r^*)$.
  \eop
  \end{enumerate}
\end{thm} 
{\bf Remarks:}   In the CoP regime, using the well-known Bendixson's criterion,  one can  prove the absence of  the limit cycles when additionally  {\small$
    \max\{u/2, \alpha\} < {\bar \beta}  < \alpha + l_i.$}  In  fact,  as in I3N case, no limit cycles are observed,  and the dynamics converge only to the LA $(i^*, r^*)$, in all our  numerical  studies \TR{(see \cite{TR})}.

\ignore{
{\small $u  (1-2i^*-c^r_i i^*) \le  \bar \beta $},  the unique local attractor is  

If further condition \eqref{Eqn_NoLC_condition} holds, then we dont' have limit cycles. The ODE solution either converges to $(i^*, r^*)$ or $(0,0)$ depending upon the initial conditions. When  \eqref{Eqn_NoLC_condition}  is not satisfied, the ODE solution can (depending upon i.c.,) have Limit cycles.

b) When {\small $u  (1-2i^*-c^r_i i^*) >  \bar\beta $}, trajectory can have following scenarios:
{\begin{enumerate}[(i)]
    \item {\bf{Global Attractor(GA)}:} In this case, trajectory $(i(t), r(t)) \to \{(i^*, r^*), (0,0)\}$ for all i.c. in $B^{(2)}.$
    \item {\bf{Local Attractor or Limit Cycle(LA-LC)}:} In this case, $(i(t), r(t)) \to \{(i^*, r^*), (0,0), \mbox{ limit cycle}\}$ for all i.c. in $B^{(2)}.$
    {\color{red} Conjucture for limit cycle}
    \item {\bf{Limit Cycle(LC)}:} In this case, trajectory $(i(t), r(t)) \to \{(0,0),\mbox{ limit cycle}\}$ for all i.c. in $B^{(2)}.$
    
\end{enumerate}}}

The possibility to eradicate the disease remains the same when  interest in the news does not change with  infection levels (see Theorems \ref{thm_sir_attractor} and   \ref{sir_linear_beta}).  Disease free state $(0,0)$ is an LA  iff $\bar\beta < \alpha$, irrespective of whether $ u  < 0$ (i.e., authentic news is predominant) or $ u  >0$ (i.e., fake news is predominant). 

\begin{figure}
\vspace{-1mm}
    \centering
    \begin{minipage}{4.1cm}
    \vspace{-5mm}
\includegraphics[width = 4.2cm, height = 4.cm]{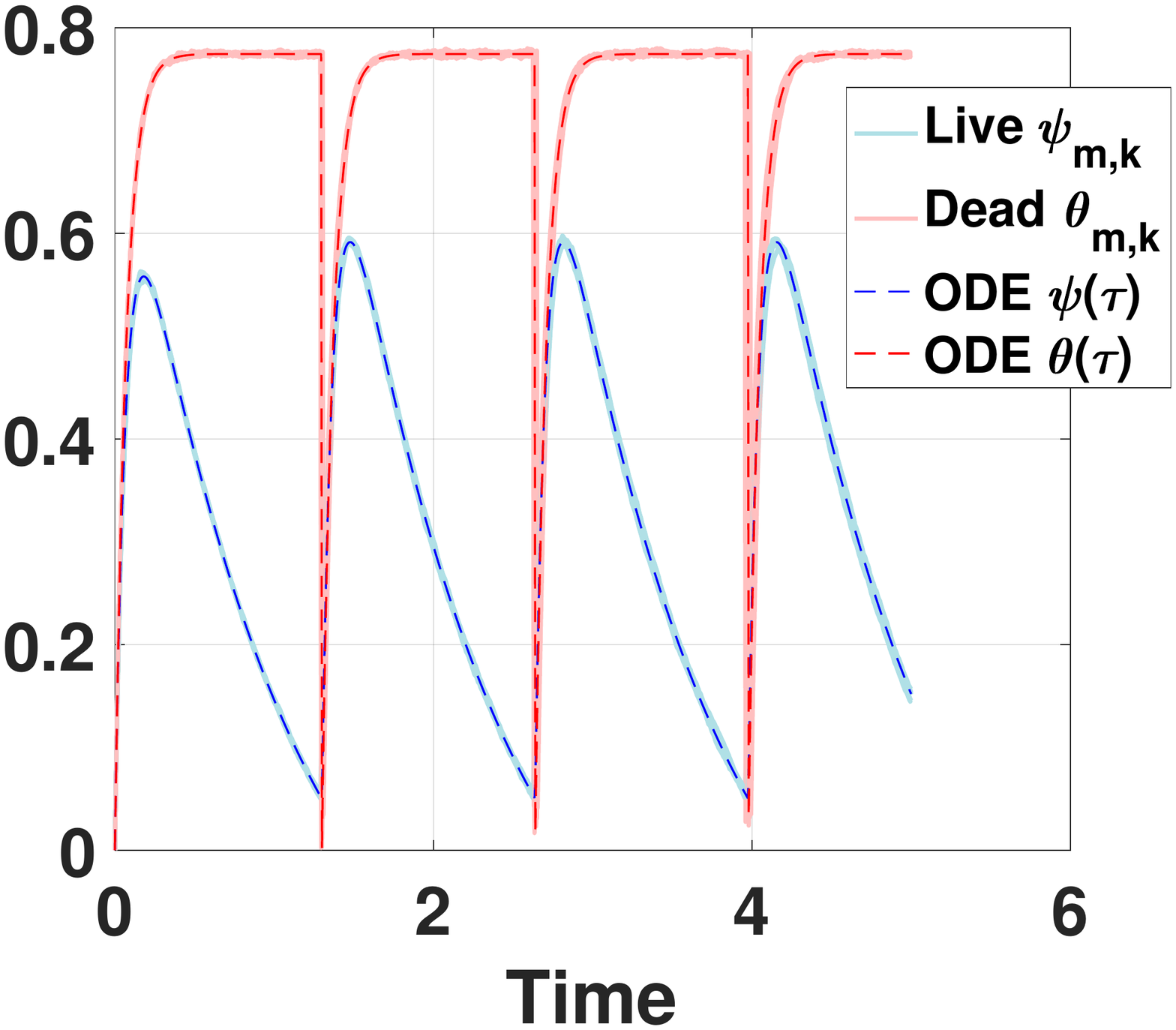}
\vspace{-14mm}
    \caption{Limit Cycle: News propagation process}
    \label{fig:limit_cycle_News}
    \end{minipage}
    \hspace{2mm}
        \begin{minipage}{4.1cm}
\includegraphics[width = 4.2cm, height = 4.cm]{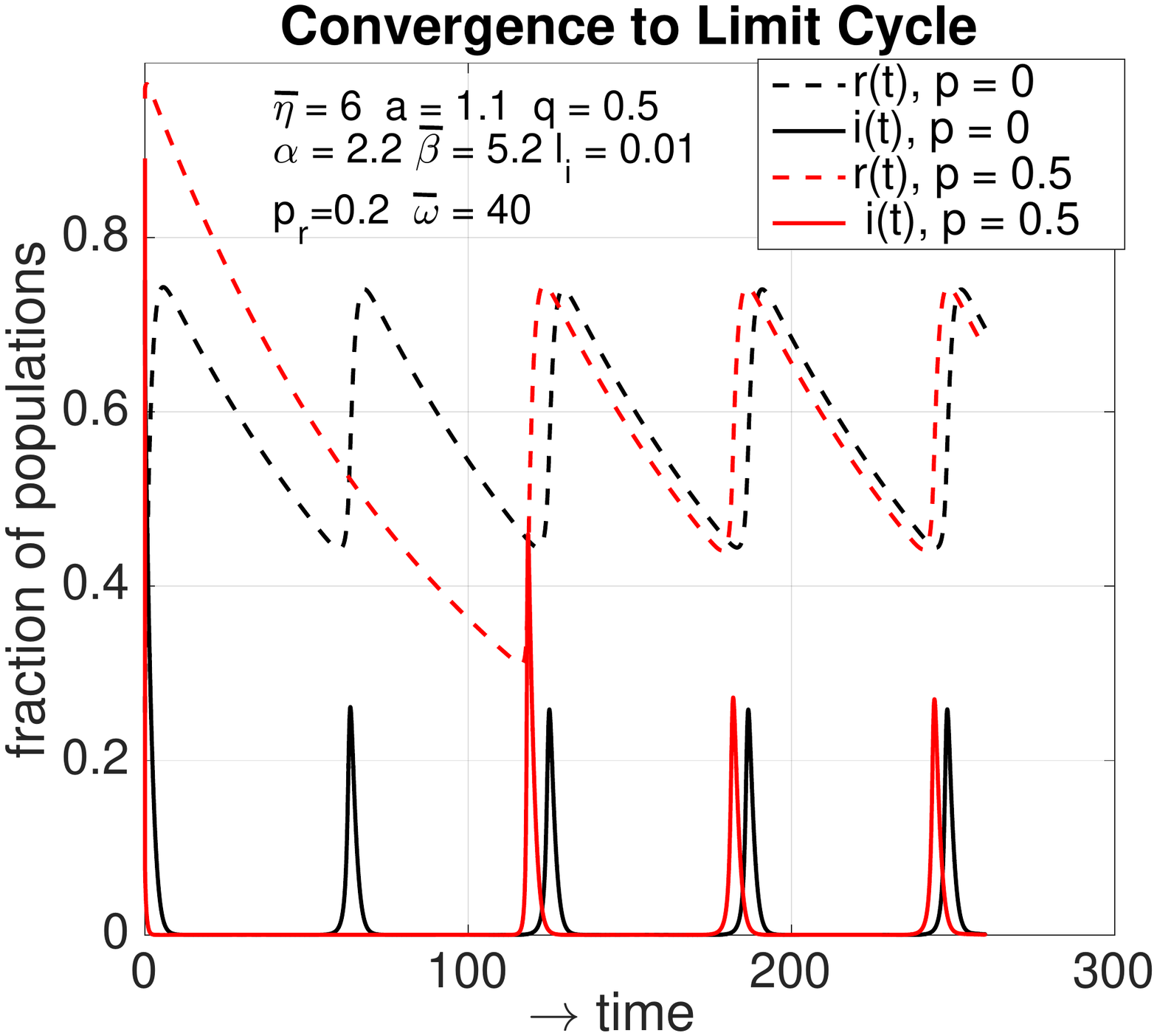}
\vspace{-14mm}
    \caption{Limit Cycle: Consolidated process. IBIN and I3N\&BIN  with $u(1-i^*-c_i^r i^* ) > \bar\beta +   l_i/{i^*}$}
    \label{fig:limit_cycle}
    \end{minipage}
    \vspace{-6mm}
\end{figure}

However, if the disease survives, the limit cycles exist when the population-activities are significantly influenced by the news (see Theorem \ref{sir_linear_beta} with  {\small $u(1-i^*-c_i^r i^* )> \bar\beta +   l_i/{i^*}$}). This result explains the waves of infections (see also Fig. \ref{fig:limit_cycle}):
a) when   people become reckless and if this reckless behavior increases with  infection level under the influence of fake news-items $(u>0)$,   the infection rises sharply to a high value and soon susceptible sub-population becomes negligible; b) then the infection starts to reduce and some recovered fraction also begins to lose immunity;  however, c) once the susceptible sub-population reaches a reasonable level, the infection rises sharply, again due to people's reckless activities induced by circulating fake news-items.

\ignore{

a) when the infection level increases, people tend to be more  careful upon receiving sufficiently authentic guidelines; b) this behavior tends to reduce the epidemic spreading;  c) however, with decrease in the infection levels, the population reverts back  to reckless behaviors, ignoring the  guidelines; and d) the reckless behaviors once again lead to an increased infection level.     
}
 
 \subsection{Increased Interest  and Behavioral Influence (IBIN$\&$I3N) }  We now consider a scenario  
in which the population seeks news more proactively and its behavior is also influenced by the news when the infection level increases.
The consolidated ODE \eqref{Eqn_combined_ODE} for this case takes the following form:
 
 \vspace{-4mm}
 {\small
 $$\frac{di}{dt} = i \left ( \left ( \bar\beta + ({\bar w}+{u} i) \frac{\bar\eta (p i + q) }{ a \bar \eta (pi+q) + 1} \right ) (1-i-r) - \alpha \right ).
 $$}
 One can have limit points as well as limit cycles as in the IBIN case. We consider one such  example  in  Fig. \ref{fig:limit_cycle}\RM{, where we further illustrate the differences between IBIN ($p=0$) and IBIN$\&$I3N ($p=0.5$) scenarios, when {\small$\bar \eta_{_{IBIN\&I3N}}=  \bar \eta_{_{IBIN }} q  =  \bar \eta$} and when ${\bar w}_{_{IBIN\&I3N}} = {\bar w} {\bar \eta}/(a {\bar \eta} + 1)  $. When the rest of the parameters are the same, \emph{both models converge to the same limit cycle.}
 In other words, the asymptotic outcome  is not dependent on whether  people consume  content with an increased interest or not, once the parameters are balanced appropriately.}  

\vspace{-3mm}
\section{Numerical Experiments}
\label{sec_numerical}

\begin{figure}[htbp]
\vspace{-3mm}
\centering 
\subfigure[Historical and fitted infection levels]{
\label{imitation figure}
\includegraphics[width=4.3cm,height = 2.3cm]{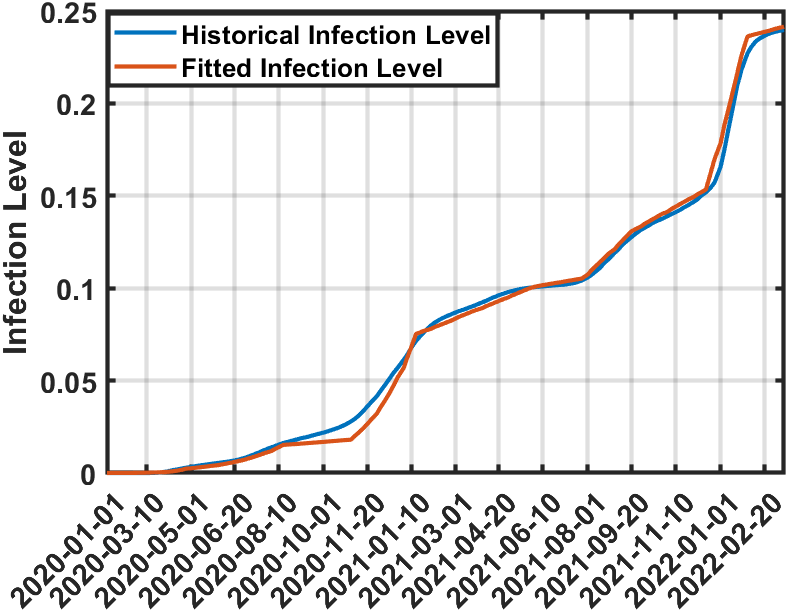}}\subfigure[Influence factors by fitting]{
\label{imitated influence factor}
\includegraphics[width=3.8cm,height = 2.3cm]{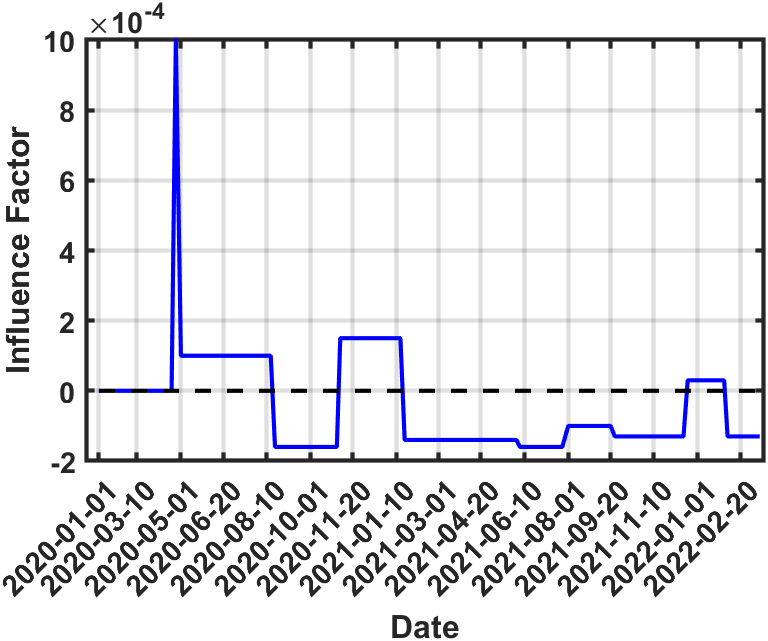}}
\vspace{-2mm}\caption{
Total cases   versus time:
Fitting  the historical COVID    datasets}
\label{ODE Approximation of News Propagation}\vspace{.2mm}
\end{figure}

In this section, we utilize the historical Twitter and infection datasets of COVID-19 to validate the relationships between the news and epidemic propagation. The real-life COVID-19 situation changes over time and is non-stationary. Building on the proposed framework, we first use the dynamics of the I3N scenario to mimic the trajectory of  the COVID-19 infection level from Jan. 1, 2020 to March 15, 2022. We inject authentic or fake news-items with different influence factors  at discrete times and obtain the trajectory of influence factors $\{\bar w\}$. Let $\bar\eta = 10$, $p = 0.7, q = 0.2,  a = 1.2$,  $\bar\beta = 0.0002, \alpha = 0.0001, i(0) = 0.00001, p_{r} = 0.2, l_{i} = 0.01.$  Fig. \ref{imitation figure} shows that we successfully imitate the historical infection levels. Fig. \ref{imitated influence factor} illustrates the influence factors of the news injected at different times. 

From Fig. \ref{imitation figure}, we observe that the COVID-19 epidemic propagates with a slow start from Jan. 15, 2020, to March 10, 2020. Then, COVID-19 bursts, and the infection level increases faster. It spawns a large amount of epidemic-related news. Due to the abrupt spreading of the epidemic and the nature of the panic-stricken population, many news-items are fake and produce a significant influence on the epidemic during the period from March 10, 2020, to May 1, 2020, in Fig. \ref{imitated influence factor}. Once the infection level increases to a notable level, public health agencies make efforts to spread authentic news, including advertising the precautionary measures and reporting the epidemic situation. Authentic news helps reduce the growth rate of the epidemic and produces a negative influence factor on the epidemic. The competition between fake and authentic news lasts for multiple cycles, shown in the time range from May 1, 2020, to Jan. 10, 2021, in Fig. \ref{imitated influence factor}. After the authentic news dominates for the period from Jan. 10, 2021, to Nov. 10, 2021, in Fig. \ref{imitated influence factor}, the public becomes accustomed to the epidemic and ignores the protective measures, rekindling the wide propagation of the epidemic and further producing more fake news in the following period from Nov. 10, 2021, to Jan. 1, 2022. The step-wise behavior of the influence depicted in Fig.  \ref{imitated influence factor} also implies a time-scale separation between the dynamics of the epidemics and the news spreading.  

\begin{figure}[htbp]
    \centering
  \vspace{-8mm}   \includegraphics[scale=0.23]{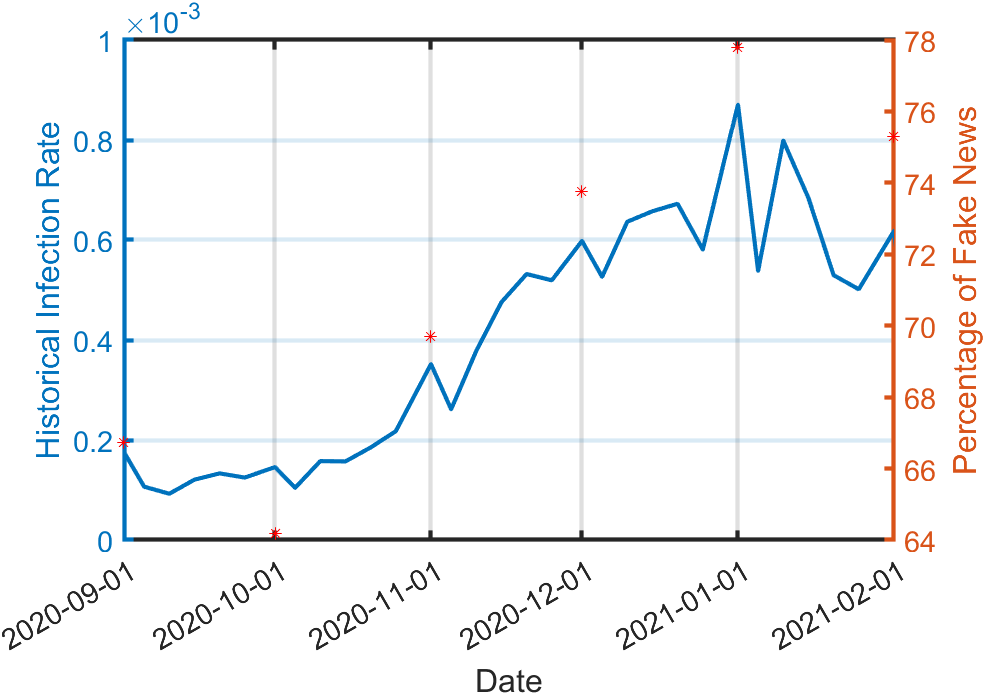}
  \vspace{-3mm}  \caption{Relationship between the historical infection rate and the percentage of fake news on Twitter. Red stars indicate the percentage of fake news on Twitter and the blue line shows the historical infection rate.}
    \label{fig:fake_news_ratio}\vspace{-1mm}
\end{figure}

We analyze the COVID-related Twitter datasets from Sept. 1, 2020, to Feb. 1, 2021, and obtain the percentage of fake news using a BERT-based approach to corroborate the results on the interplay between the epidemic and fake news. 
Fig. \ref{fig:fake_news_ratio} shows that the percentage of fake news on Twitter is positively correlated with the historical infection rate. \RM{The red stars indicate the percentage of fake news on Twitter at discrete time points. The blue line shows the historical infection rate. When the epidemic becomes severer, it generates public panic and promotes the spreading of fake news. Hence, the percentage of fake news on Twitter increases. However, when the epidemic becomes less aggressive, the percentage of fake news decreases.}

The understanding of the integrated dynamics enables a short-range prediction despite the non-stationarity in the long run.

\vspace{-5mm}
\section{Conclusions}\label{sec:Conclusions}

In this work, we have presented a two-timescale integrative epidemic model that consolidates the propagation of online social network (OSN) news-items into the compartmental epidemic models; the aim is to study the impact of fake news on the spread of infectious disease. We have analyzed the equilibrium behavior of the coupled dynamic system under different behavioral responses to fake news. We found a periodic pattern of the epidemic when the behavioral response is linearly dependent on the infection level. Our model explains multiple waves of infection and reinfection as seen in the COVID-19 pandemic. We use Twitter and COVID-19 datasets to validate the proposed model.   Posting appropriate authentic information can help eradicate the severe spreading of the disease. Misinformation can prolong the epidemic. The successful fitting of the proposed model with the historic COVID-19 dataset shows that the incorporation of human behaviors into the epidemic model provides a promising approach to predicting the trends of the epidemic processes and evaluating public health policies. Motivated by this observation, our future work would focus on the optimal design of intervention mechanisms and policies to mitigate the impact of the pandemic.    

 
 \RM{
 We have witnessed many pandemic diseases in the past (including COVID-19), and they exhibit natural period behavior due to the time-varying characteristics of the new mutants. However, the existence of OSNs which facilitate the rapid transfer of information across the globe in the current pandemic has a different impact and the goal of the paper is to capture this   aspect. More precisely, the paper proves the existence of infections and reinfections even without the changes in the characteristics of the mutants, which are predominately influenced by the circulating information. 
 }
 
\vspace{-2mm}

\bibliographystyle{IEEEtran}
\bibliography{sample}

\RM{

\vspace{-3mm}
\section*{Appendix}
\vspace{-2mm}

\TR{}
{
\noindent {\bf Proof of Theorem \ref{newsode}:} 
We consider the general system, \eqref{Eqn_update_stocahstic}, where assumptions {\bf A.1} and {\bf A.2} are satisfied. We show that the equations given by \eqref{Eqn_update_stocahstic} can be approximated by the solution of the ODE \eqref{Eqn_update_ODE}.
We prove this result
  using  \cite[Theorem 9, pp. 232]{benven}.

  {\bf Assumption required by \cite[Theorem 9, pp. 232]{benven} }
  To this end, we need to show that the relevant assumptions are satisfied, which are reproduced below as {\bf B}.1-{\bf B}.4, in our notations.  Consider a stochastic iterative scheme of the type:
  \begin{eqnarray}
  \label{Eqn_SA_scheme}
  \phi_{k+1} = \phi_k + \epsilon 
  H (\xi_{k+1} , \phi_k ).
  \end{eqnarray}
If the above scheme satisfies {\bf B}.1-{\bf B}.4, then \cite[Theorem 9, pp. 232]{benven} is applicable.

\begin{enumerate}[{\bf B}.1]
    \item 
 There exists a family ${P_\phi}$ of transition probabilities ${P_\phi (\xi, A)}$ such that, for any Borel subset $A$, we have
\begin{eqnarray}
P[\xi_
{n+1} \in A \mid \mathcal{F}_n] = P_{\phi_n}(\xi_n, A),
\end{eqnarray}
where $\mathcal{F}_n = \sigma(\phi_0, \xi_0, \xi_1, \dots , \xi_n)$ is a sigma-algebra. This in turn implies that
the tuple $(\xi_k, \phi_k)$ forms a Markov chain.\newline

\item  For any compact subset $Q$ of $D$, there exist constants $C_1, q_1$ such that for all $\phi \in Q $  we have
\begin{eqnarray}
\|H (\xi, \phi)\| \le C_1 (1+ | \xi|^{q_1}).
\end{eqnarray}

\item  There exists a function $h_1$ on $D$, and for each $\phi \in D$ a function $\nu_\phi (\cdot)$ such that
\begin{enumerate}
    \item $h_1$  is locally Lipschitz on $D$.
    \item $(I - P_{\phi}) \nu_{\phi} (\xi) = H (\xi, \phi) - h_1 (\phi)$, where $I$ is the identity matrix of the same order as the one of $P_{\phi}$.
    \item For all compact subsets $Q$ of $D$, there exist constants $C_2, C_3, q_2, q_3$ and $\lambda \in [0.5,1]$, such that for all $\phi , \phi' \in Q$, following is true:

\begin{enumerate} [i.]
\small
    \item $\| \nu_\phi (\xi)\| \le C_2 (1+ | \xi |^{q_2})$,
    
    \item $\| P_\phi \nu_\phi (\xi) - P_\phi' \nu_\phi' (\xi)\| \le C_3 (1+ | \xi |^{q_3}) \| \phi - \phi'\|^\lambda.$
\end{enumerate}
\end{enumerate}

\item  For any compact set $Q$ in $D$ and for any $q>0$, there exists a $\mu_q(Q)<\infty$, such that
for all n and $\xi,  \phi \in R^d$, following is true:

\vspace{-4mm}

{\small
\begin{eqnarray}
\hspace{-7mm} E_{\xi, \phi} \{ I(\phi_k , \xi_k \in Q, k \le n) (1 + \|\xi_{n+1}\|^q)\} \nonumber   
\\ &\hspace{-40mm}\le  \mu_q (Q) (1 + \|\xi_{n}\|^q), 
\label{Eqn_exp}
\end{eqnarray}}
where $E_{\xi, \phi}$ represents the expectation taken with initial conditions $(\xi_0 ,\phi_0) = (\xi , \phi)$.
\end{enumerate}

{\bf Assumptions :}  We now prove that the above assumptions are satisfied by  \eqref{Eqn_update_stocahstic}. 
First observe that \eqref{Eqn_update_stocahstic} has the same form as in 
\eqref{Eqn_SA_scheme}, with  $\phi = [\psi, \theta]^T$ and,
 \begin{eqnarray}
 \label{Eqn_H_xi_phi}
H (\xi, \phi) = \begin{bmatrix}
\mathcal{J}_{\phi} ( \xi   - \theta) + (1 - \mathcal{J}_{\phi})(-C\theta)\\
\mathcal{J}_{\phi} ( \xi   - \psi - 1) + (1 - \mathcal{J}_{\phi})(-\psi)
\end{bmatrix}.
\end{eqnarray}

We now prove the required assumptions one after the other.  
\begin{enumerate}[(i)] 
\item The offsprings $\xi_{k+1}$ depend only upon the total population $\theta_k$ and hence
 assumption {\bf B}.1 is satisfied with
\begin{eqnarray*}
P_{\phi_n} (\xi_n, A) &=& Q_{\theta_n} (A) \mbox{, where,} \\ 
\  Q_{\theta_n} (A) &:=& P ( \xi_{n+1} \in A | \theta_n ).
\end{eqnarray*}

\item  It is direct from \eqref{Eqn_H_xi_phi}
that for any compact $Q$
\begin{eqnarray}
\|H (\xi, \phi)\| \le \sup_{\phi \in Q} \| \phi \| (1+ |\xi|)
\end{eqnarray}
{\bf B}.2 is satisfied with
  $q_1 = 1$ and $C_1 = \sup_{\phi \in Q} ||\phi||$ .

     \item We will show that assumption {\bf B}.3 is satisfying by setting $\nu_\phi (\xi) := H (\xi, \phi)$ and 
\begin{eqnarray*}
     h_1 (\phi) &:=& (P_{\phi})\nu_{\phi}  (\xi)  \mbox{, where } \\
(P_{\phi})\nu_{\phi}  (\xi) &=&  \int H (y, \phi)(Q_{\phi})dy,
\end{eqnarray*} Observe that under {\bf A}.1-2, we have: 

\begin{eqnarray}
h_1 (\phi) = \begin{bmatrix}
\mathcal{J}_{\phi} (\mathcal{M}(\theta) - \theta) + (1 - \mathcal{J}_{\phi})(-C\theta)\\
\mathcal{J}_{\phi} (\mathcal{M}(\theta) -1 - \psi) + (1 - \mathcal{J}_{\phi})(-\psi)
\end{bmatrix}.
\label{h_1}
\end{eqnarray}
We will now prove all the sub-assumptions {\bf B}.3.a-c in the following:
\begin{enumerate} [{\bf a.}]
    \item From equation \eqref{h_1}, under \textbf{A}.2, we have:
     {\small
    \begin{eqnarray} 
    \label{h_1_lips}
     | h_1 (\phi) - h_1(\phi')| \hspace{-6mm} \nonumber &
  \\ & \hspace{-10mm}\le (C + 1) |\theta - \theta'| + 2|\psi - \psi'| + |\mathcal{M}(\theta) - \mathcal{M}(\theta')| \nonumber\\
  &\hspace{-29mm}\stackrel{a}{\le} (C + 1 + 2\rho) |\theta - \theta'| + 2|\psi - \psi'|\\
  &\hspace{-1.6cm}\le D |\phi - \phi'|, \mbox{ with } D := \max \{ C + 1 + 2\rho, 2 \} \nonumber,
    \end{eqnarray}
   }where `a' hold true as $\lambda \leq 1$. Thus,  $h_1$ is locally Lipschitz.
    \item  The definitions of $h_1 (\phi)$ and $H(\xi, \phi)$ make this obvious.
    \item 
    \begin{enumerate}[i.)]
        \item  This proved along with {\bf B}.2,
   as
   \begin{eqnarray*}
   \nu_\phi(\xi) = H (\xi, \phi).
   \end{eqnarray*}  

   \item From \eqref{h_1_lips} and definition of $h_1(\phi)$, this assumption is satisfied.

  \end{enumerate}
\end{enumerate}

\item For proving {\bf B}.4, consider any compact $Q$, then we can upper bound the LHS of \eqref{Eqn_exp} as below:
\begin{align*}
 E_{\xi, \phi} \{ I(\phi_k , \xi_k \in Q, k \le n) (1 + \|\xi_{n+1}\|^q)\} \\ &\hspace{-40mm}\le E_{\xi, \phi} \{I(\phi_n , \xi_n \in Q) (1 + \|\xi_{n+1}\|^q)\}\\
 &\hspace{-40mm}\le \sup_{\phi \in Q}  \int (1 + y^q) Q(\phi,dy)\\ 
 &\hspace{-40mm}= \sup_{\phi \in Q} E[1+\xi_{n+1}^q | \theta_n = \theta]\\ 
&\hspace{-40mm} \stackrel{a}{\le} 1 + \rho^q \\
&\hspace{-40mm} \le (1 + \rho^q) (1+ ||\xi_n ||^q)
\end{align*}
for any $q > 0$, 
where  equality `a' follows from assumption {\bf A}.2. Then, \textbf{B}.4 hold with $\mu_q(Q) := 1 + \rho^q$.

\end{enumerate} 
Now, using  \cite[Theorem 9, pp. 232]{benven}, the result is proved. 
 \eop
}

\noindent {\bf ODE Solution and the Derivation of \eqref{limit_of_frac}:} 
We consider the system \eqref{Eqn_update_stocahstic} and  the corresponding ODEs are given in \eqref{Eqn_update_ODE}.

Say, we start in {\bf regime $\mathcal{R}^1$}, and say  with $\psi(0) = \delta_\psi$ and some $\theta(0) > \delta_\theta$. Then, the solution is given by  

\vspace{-5mm}
{\small \begin{eqnarray*}
\psi(\tau) = \delta_\psi e^{-\tau} \mbox{ \normalsize and }  \theta(\tau) = \theta(0) e^{-C\tau}.
\end{eqnarray*}}
And then, there exists a time epoch (say $\tau_1 = \nu_1$) such that
{\small
\begin{eqnarray*}
\theta(\nu_1) = \delta_\theta \  \  (\mbox{\normalsize say } \theta_{01})  \mbox{ \normalsize and let } \psi(\nu_1) := \psi_{01}.
\end{eqnarray*}
}
{\bf Regime ${\cal R}^2$:} Let $\Delta = a\eta + 1$. 
The solution after
   time epoch $ \tau_1 $, 
with initial conditions $\phi(\tau_1) = (\psi_{01}, \delta_\theta)$,
is given by:

\vspace{-4mm}
{\small
\begin{align*}
\psi(\tau +\tau_1) &=  \Delta \theta(\tau + \tau_1)    +(\psi_{01} - \delta_\theta \Delta  + a \eta^2 +1) e^{- \tau}  -  a \eta^2 -1, \nonumber \\
 \theta (\tau + \tau_1) &= \left(  \delta_\theta - \frac{\eta}{\Delta}   \right) e^{-\Delta \tau} + \frac{\eta}{\Delta }.  
 \label{Eqn_aftertime1}
\end{align*}
}Then, there exists  epoch {\small $\tau_2 = \nu_2 + \tau_1 $}, with {\small $ \nu_2 > 0 $}, such that

\vspace{-3mm} 
{\small
\begin{eqnarray*}
 \psi(\tau_2) &=& \delta_\psi =: \  \psi_{02} \mbox{ \normalsize and } \\
 \theta(\tau_2) &=& \left(  \delta_\theta - \frac{\eta}{\Delta }   \right) e^{-\Delta \nu_2} + \frac{\eta}{\Delta } 
 =: \theta_{02}.
\end{eqnarray*}
}This completes $1$-st cycle, i.e., $L_{C_1}$. Now it goes back to regime 1, with initial conditions $\phi(\tau_2) = (\psi_{02}, \theta_{02} )$ and continues this pattern. We can construct series of time epochs $\tau_i$ with appropriate initial conditions, constructed from the terminal conditions of previous regimes. 

This process continues and we arrive at a periodic asymptotic solution. This limit can be specified if the initial/terminal conditions can be identified by solving the following fixed-point equations (recall $\Delta  = a \eta +1 $):

\vspace{-4mm}
{\small 
\begin{eqnarray*}
\theta_{02}^* = \theta (\nu_2^*) &=& \left(  \delta_\theta - \frac{\eta}{\Delta }   \right) e^{-\Delta \nu_2^*} + \frac{\eta}{\Delta }, \\
\psi_{01}^* = \psi (\nu_1^*) &=& \delta_\psi e^{-\nu_1^*} \mbox{, where, }\
\theta_{02}^* e^{-C  \nu_1^*} \  = \  \delta_\theta \mbox{ and } \\
&& \hspace{-30mm}   \Delta \theta_{02}^*   
 +(\psi^*_{01} - \delta_\theta \Delta + a \eta^2 +1) e^{- \nu_2^*} - a \eta^2 -1   = \delta_\psi. \nonumber 
\end{eqnarray*}}

{\bf At the limit and in Regime 2}, we have for any $\tau \le \nu_2^*$
{\small
\begin{eqnarray*}
\theta (\tau + \tau_\infty)   &=& \left(  \delta_\theta - \frac{\eta}{\Delta }   \right) e^{-\Delta \tau} + \frac{\eta}{\Delta }  .
\end{eqnarray*}
}
 When $\nu_2^*$ is sufficiently large, one can approximate $\theta_{02}^* \approx \eta/\Delta.$  Under this approximation, the fixed-point equations are solved by:

\vspace{-5mm}
{\small \begin{eqnarray*}
\theta_{02}^* \approx  \frac{\eta}{\Delta}, \  \   \ 
\nu_1^* = \frac{-1}{C} \ln \left( \frac{ \Delta \delta_\theta}{\eta} \right), 
\ 
\psi_{01}^*  = \delta_\psi e^{-\nu_1^*}
\approx \delta_\psi \left ( \frac{\Delta \delta_\theta}{\eta} \right )^{1/C}.
\end{eqnarray*}
Thus, the limit in \eqref{limit_of_frac} equals
\begin{eqnarray}
\label{Eqn_limit_of_frac_proved}
  {\theta}_\infty ^{*} (\eta) := \lim_{k \to \infty} \max_{\tau' \in L_{C_k} } \theta (\tau') = \theta_{02}^* \approx  \frac{\eta } { a \eta  + 1}.
\end{eqnarray}}

\ignore{
Similarly, at $\tau_3$, we have
\begin{eqnarray*}
\psi(\tau_3) = \delta_\psi e^{-\tau^*}, \mbox{ where } e^{-C\tau^*} = \frac{\delta_\theta (a \eta +1)}{\eta}.
\end{eqnarray*}
On continuing like this, we obtain the result at the limit, and we will switch firstly from Phase 1 at a time epoch (say at $\bar\tau'$) such that
\begin{eqnarray*}
\psi (\bar\tau') = \delta_\psi \mbox{ and } \theta (\bar\tau') = \frac{\eta}{a \eta +1 }.
\end{eqnarray*}
Then, we will switch after some time (say $\tau_1^*$) such that 
\begin{eqnarray*}
\theta(\bar\tau' + \tau_1^*) = \delta_\theta, \mbox{ which implies } \tau_1^* = \frac{-ln \left( \frac{(a\eta +1 )\delta_\theta}{\eta} \right) }{C}.
\end{eqnarray*}
So, Phase 1 ends with length $\tau_1^*$ such that 
\begin{eqnarray}
\hspace{-15mm}\psi(\tau_1^*) = \delta_\psi e^{-\tau_1^*} \mbox{ and }   \theta(\tau_1^*) = \delta_\theta. \nonumber
\end{eqnarray}
  Now from \eqref{Eqn_aftertime1}, Phase 2 ends after time $\tau_2^*$ such that 
  \vspace{-1mm}
  {\small
 \begin{eqnarray}
 \theta (\tau_2^* + \tau_1^*) &=& \left(  \delta_\theta - \frac{\eta}{a \eta +1 }   \right) e^{-(a \eta +1)\tau_2^*} + \frac{\eta}{a \eta +1 } \mbox{ \normalsize  and } \nonumber \\ \psi (\tau_2^* + \tau_1^*) &=&\delta_\psi \nonumber
 \end{eqnarray}
  }    

and $\tau_2^*$ can be computed by solving the following fixed-point equation
\begin{eqnarray*}
\theta (\tau_2^*) &=& \left(  \delta_\theta - \frac{\eta}{a \eta +1 }   \right) e^{-(a \eta +1)\tau_2^*} + \frac{\eta}{a \eta +1 }.
\end{eqnarray*}
}

\section*{Appendix:  Proofs related to Epidemic-ODE   }
It is easy to verify that global unique solutions exist for all ODEs considered, as the Right-Hand Sides (RHS) are Lipschitz continuous.
We first show that  for any initial condition (i.c.) in $\B$, the solution remains in ${B}^{(2)}$:

\begin{lem}
\label{lemma1}
Consider the two-dimensional ODE  defined in \eqref{Eqn_combined_ODE}.
\ignore{
{\small \begin{eqnarray*}
\frac{di}{dt} &=& i \left ( \left ( \bar\beta + ({\bar {  w}}+{{ u}} i) \frac{\bar\eta (p i + q) }{ a \bar \eta (pi+q) + 1} \right ) (1-i-r) - \alpha \right )\\ 
\frac{d r}{d t} &=&  \left  (   i \alpha p_r - r l_i\right )
\end{eqnarray*}}} 
Consider any  i.c. $(i_0, r_0) \in { B}^{(2)} $, 
 then the unique solution $(i(t), r(t)) \in B^{(2)}$ for all $t$.
\end{lem} 
{\bf{Proof:}}
\nchanges{
We will show that  
the flow of the ODE on the boundary of $\B$ is either inwards or on the boundary (see Fig. \ref{fig_lemma1}).  This provides the required result by Nagumo’s theorem \cite{ghorbal2022characterizing}, which then says that  any trajectory started inside region $\B$ will remain in $\B$. 

\begin{figure}[h]
    \centering
  \vspace{-5mm}   \includegraphics[scale=0.11]{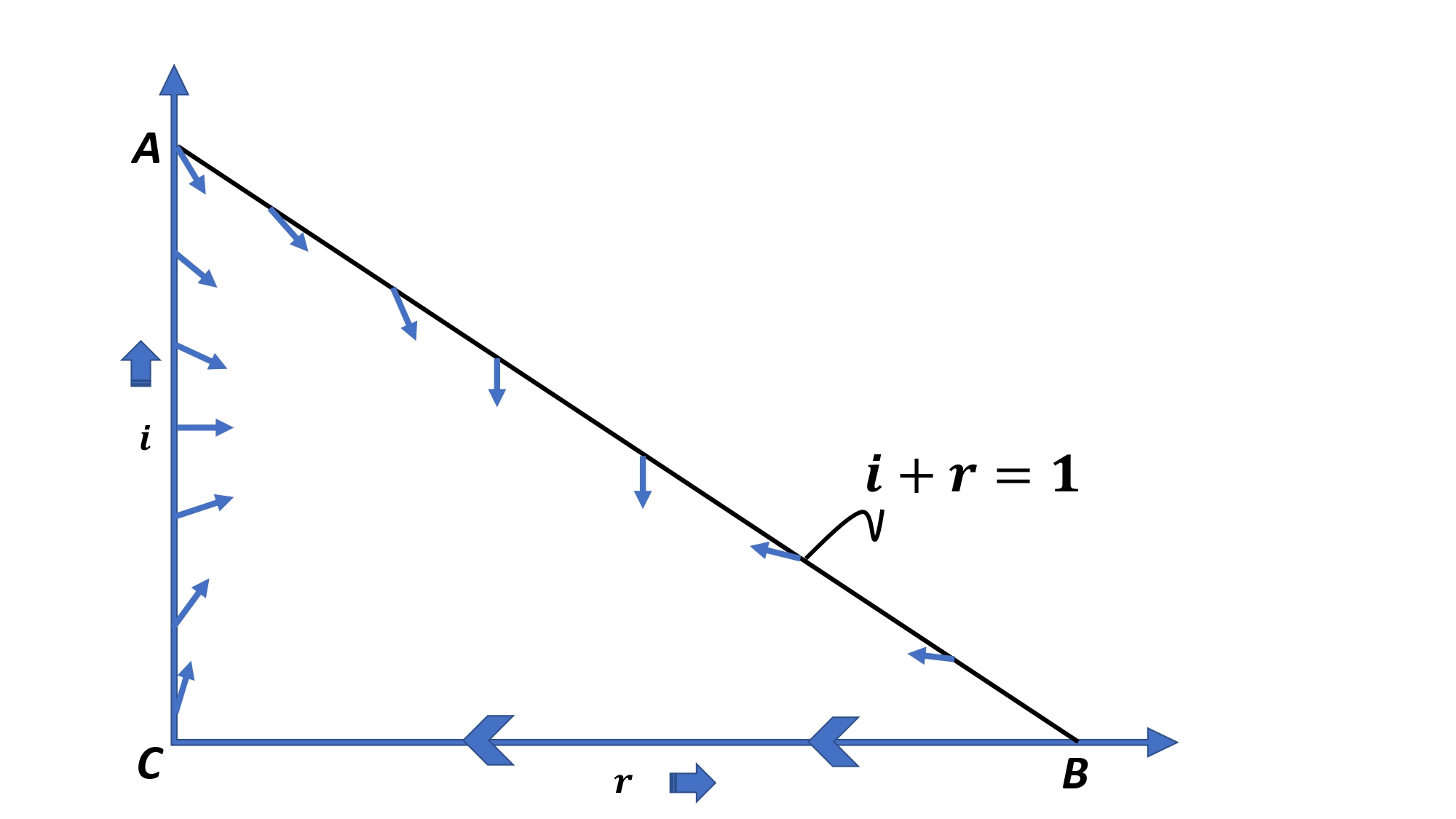}
  \vspace{-3mm}  \caption{ODE Trajectory}
    \label{fig_lemma1}\vspace{-4mm}
\end{figure}

Consider any  $(r,i) \in \partial \B$. When the boundary point is $(r, 0)$ (i.e., on $r$-axis as shown in  figure),  then clearly from ODE \eqref{Eqn_combined_ODE}, $di/dt = 0$ while    $dr/dt < 0$, so the flow is directed towards $(0,0)$ and along the boundary line $BC$. 
Now, for any point $(0,i)$ on the boundary $AC$, again from ODE \eqref{Eqn_combined_ODE}, $dr/dt > 0$ and   depending upon $\beta(i)$, the ODE  flow is either upwards and towards right or downwards and towards right  (as shown in figure \ref{fig_lemma1}); in either case the flow is into $\B$.  At the corner point $A$, the slope of flow is 
$
\frac{di}{dr}  = \frac{- \alpha i}{i \alpha p_r} = \frac{-1}{p_r},
$
whose magnitude is greater than one (because $p_r$ is the  fraction of recovered that get  immunized, and hence $p_r < 1$), which is hence ${di}/{dr}$ is smaller than the slope of $AB$ (as the slope of line $i+r =1$ is $-1$). So, flow at $A$ is also into $\B$. Now consider any point $(r, i) $ on line $AB$; here $i+r = 1$ and hence the slope of the flow,
$
\frac{di}{dr} = \frac{- \alpha i}{i (\alpha p_r +l_i) -  l_i },
$
which is increasing   with $i$ and hence ${di}/{dr}  $ for any $(r,i)$ on $AB$
is less than that at $A$, i.e., ${di}/{dr} < -1/p_r$.
Thus the flow is   inward on $AB$ also.  Proof follows  by Nagumo’s theorem \cite{ghorbal2022characterizing}. \eop
}


\noindent {\bf Proof of Theorem \ref{thm_sir_attractor}:}
  For the ODE \eqref{epi_ode}, (0,0) is an equilibrium point, we have another   equilibrium point,  $(i^*,r^*)$, 
  
   \vspace{-4mm}
  { \small{\begin{eqnarray}
i^* =  \frac{ \bar\beta - \alpha }{\bar\beta(1+\frac{\alpha p_r}{l_i}
)}  \ \mbox{ and }   r^* \ = \ \frac{ i^* \alpha p_r}{l_i}. 
\label{Eqn_istar_rstar}
\end{eqnarray}}}
Furthermore, the Jacobian matrix at any $(i,r)$ is:

\vspace{-2mm}

{\small $$J (i, r) = \begin{bmatrix}
 \bar\beta  (1-i-r) - \alpha  - \bar\beta i & -i \bar\beta   \\
\alpha p_r & -l_i 
\end{bmatrix}.$$}Matrix  $J(0,0)$ is negative definite and (0,0) is   LA iff $\bar \beta <   \alpha$.
\ignore{
But when $(i,r) \to (0,0)$, then
\begin{eqnarray*}
    \frac{di}{dt} \ \to 0 ,   
    \  \  \   \ 
\frac{d r}{d t} \ \to 0.  
\end{eqnarray*}}

{\bf Case with $\bar \beta <   \alpha$:} 
For this case, from \eqref{Eqn_istar_rstar}, $i^* <0$ for the non-zero equilibrium point. Thus there is only one   equilibrium point in ${ B}^{(2)} $, which is $(0,0)$.
Further, there cannot be a limit cycle inside $\B$, as it has to enclose $(0,0)$ (and not touch it), see  \cite{palis2012geometric}, and, as any trajectory (with i.c., in $\B$) is trapped inside $\B$ by Lemma \ref{lemma1}. Thus by Poincar\'e–Bendixson theorem (e.g., \cite{palis2012geometric}),    the limit set of any trajectory staring in $\B$ is $\{(0,0)\}$, which implies that $(0,0)$ is a global attractor  (as it is a local attractor).    

\ignore{
, and we will show that the ODE trajectory $(i(t), r(t) ) \to (0,0)$ for any i.c. in ${ B}^{(2)} $. 

Observe that the derivative  of  $i(\cdot)$ is negative for this case:
$$ \frac{di}{dt} \ =  \ i \big ( \bar \beta(1-i - r)-  \alpha \big ) < \bar \beta - \alpha. $$ 
Thus  (e.g., by  \cite[Lemma 6]{agarwal2021new}),   $ i(t) \to 0  $ as $t\to \infty$ for any i.c. In fact, $i(t) \downarrow 0$. Thus,  for  any $\tau_1 > 0$ the derivative of the second component can be upper bounded as below:
\begin{eqnarray}
\label{dr}
\frac{dr}{dt} \leq i(\tau_1) \alpha p_r - r l_i, \ \forall t \geq \tau_1.
\end{eqnarray}
 Then (e.g.,by \cite{piccinini2012ordinary}), we have the following for any $t \ge \tau_1$:
 
{\small \begin{eqnarray*}
 0 \le r(t) \leq  \bigg(r(\tau_1) - i(\tau_1) \alpha p_r \bigg ) e^{-l_i (t-\tau_1)} + i(\tau_1) \alpha p_r , \ \forall t \geq \tau_1.
\end{eqnarray*}}
\ignore{
Then by comparison of ODEs, for every $\epsilon > 0 $ there exists $\tau_2 \geq 0$ such that 
\begin{eqnarray}
\label{rt}
r(t) \leq i(\tau_1) \alpha p_r + \frac{\epsilon}{2}, \ \forall t \geq \tau_2,
\end{eqnarray}}
Thus and further because $i(\tau_1) $ is decreasing as $\tau_1$ increases and because $r(\tau_1) $ is bounded uniformly by one we have that $r(t) \to 0$ for any i.c. in $B^{(2)}$.
}


\noindent {\bf Case $\bar\beta > \alpha $:} 
\ignore {In this case, (0,0) is not an attractor as one can show that the Jacobian is not negative definite at (0,0). For non zero $(i^*,r^*)$ given in \eqref{Eqn_istar_rstar}, clearly 
if we assume 
\begin{eqnarray}
\label{xy}
x(t) = i(t) - i^* \mbox{ and } y(t) = r(t) - r^*.
\end{eqnarray}
Then,
\begin{eqnarray*}
    \frac{dx}{dt} \ = \  -\bar\beta (x+i^*) (x+y) \to 0, \mbox{ and }\\   
    \  \ \ \ \ \ 
\frac{d y}{d t} \ = \  x\alpha p_r - y l_i ,\\ \mbox{ by doing similar steps as \eqref{dr}. }  
\end{eqnarray*}
This completes the proof.}
The Jacobian matrix, 
{\small 
$
J (i^*, r^*) = \begin{bmatrix}
  - \bar\beta i^* & - \bar\beta i^*    \\
\alpha p_r & -l_i,
\end{bmatrix},$}
 is negative definite  (use minors). Thus $(i^*, r^*)$ is an LA.

In Lemma  \ref{lemma2}, we  have constructed a bounded region  (e.g.,   ABCD in left sub-figure,  Fig. \ref{fig1}), call it ${\cal R}$,  with  flow inwards or on boundary. On applying Nagumo’s theorem \cite{ghorbal2022characterizing}, any trajectory started inside region ${\cal R}$ will remain in ${\cal R}$. Further, using the Bendixson criterion (e.g., \cite{verhulst2006nonlinear}), there will be no limit cycle inside the region ${\cal R}$ 
because of  \eqref{Eqn_fg} in  Lemma~\ref{lemma2}.
 
Then, by Poincar\'e–Bendixson (PB) theorem (e.g., \cite[Theorem 1.8]{palis2012geometric}), the ODE trajectory starting in ${\cal R}$  has a unique limit point,  and as the limit point is an LA, the trajectory  converges to $(i^*, r^*)$. Now, by the critical point criterion (\cite{palis2012geometric}), any closed orbit has a critical point (which is not a saddle point) in its interior.
\Attempt{
\begin{figure}[htbp]
    \centering
  \vspace{-2mm}   \includegraphics[scale=0.2]{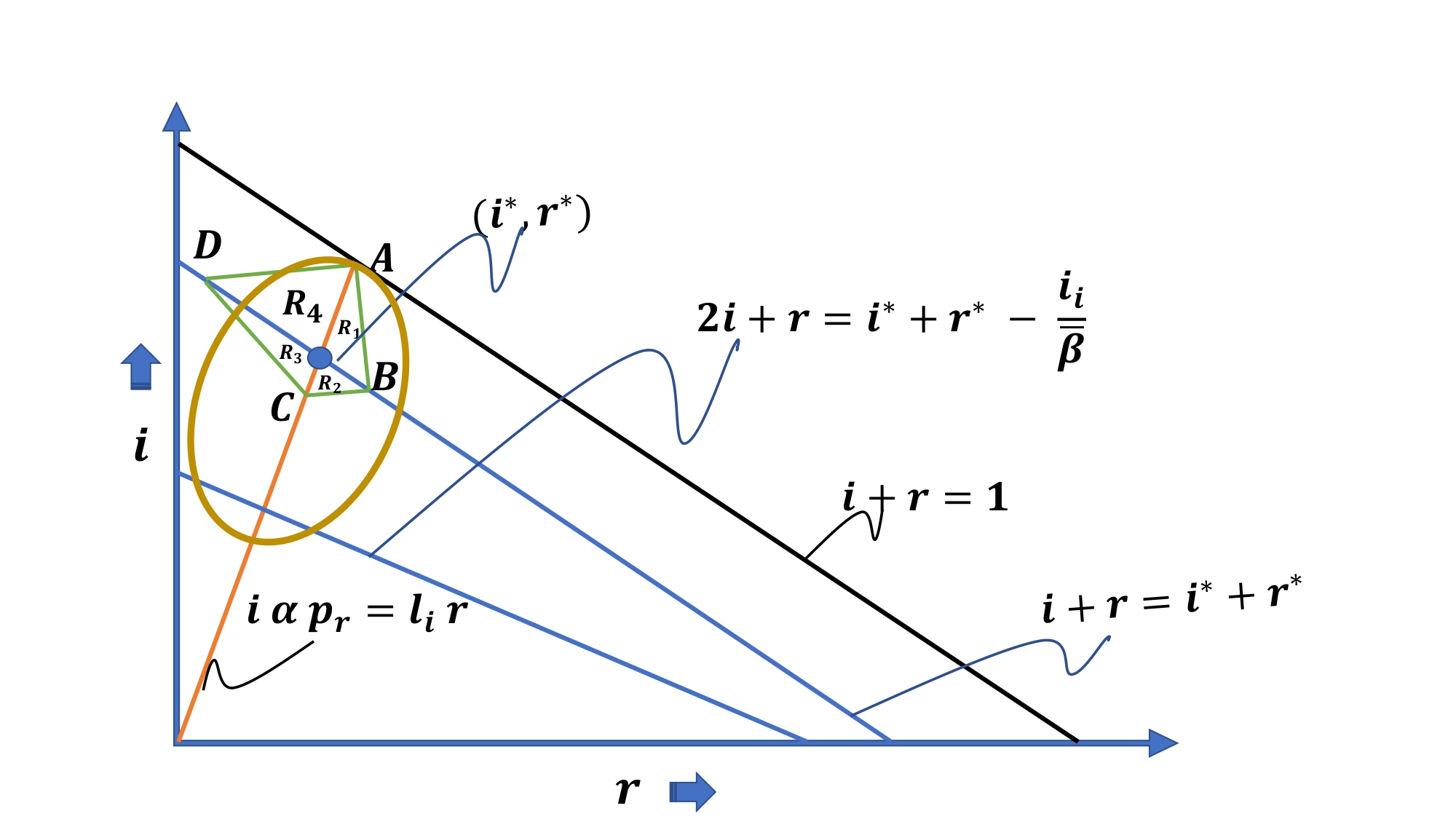}
  \vspace{-3mm}  \caption{Limit cycle when $\alpha + l_i < \bar \beta$,  $\bar \beta < 1 + \frac{l_i}{\alpha p_r}$}
   \vspace{-1mm}
\end{figure}

\begin{figure}[htbp]
    \centering
  \vspace{-2mm}   \includegraphics[scale=0.2]{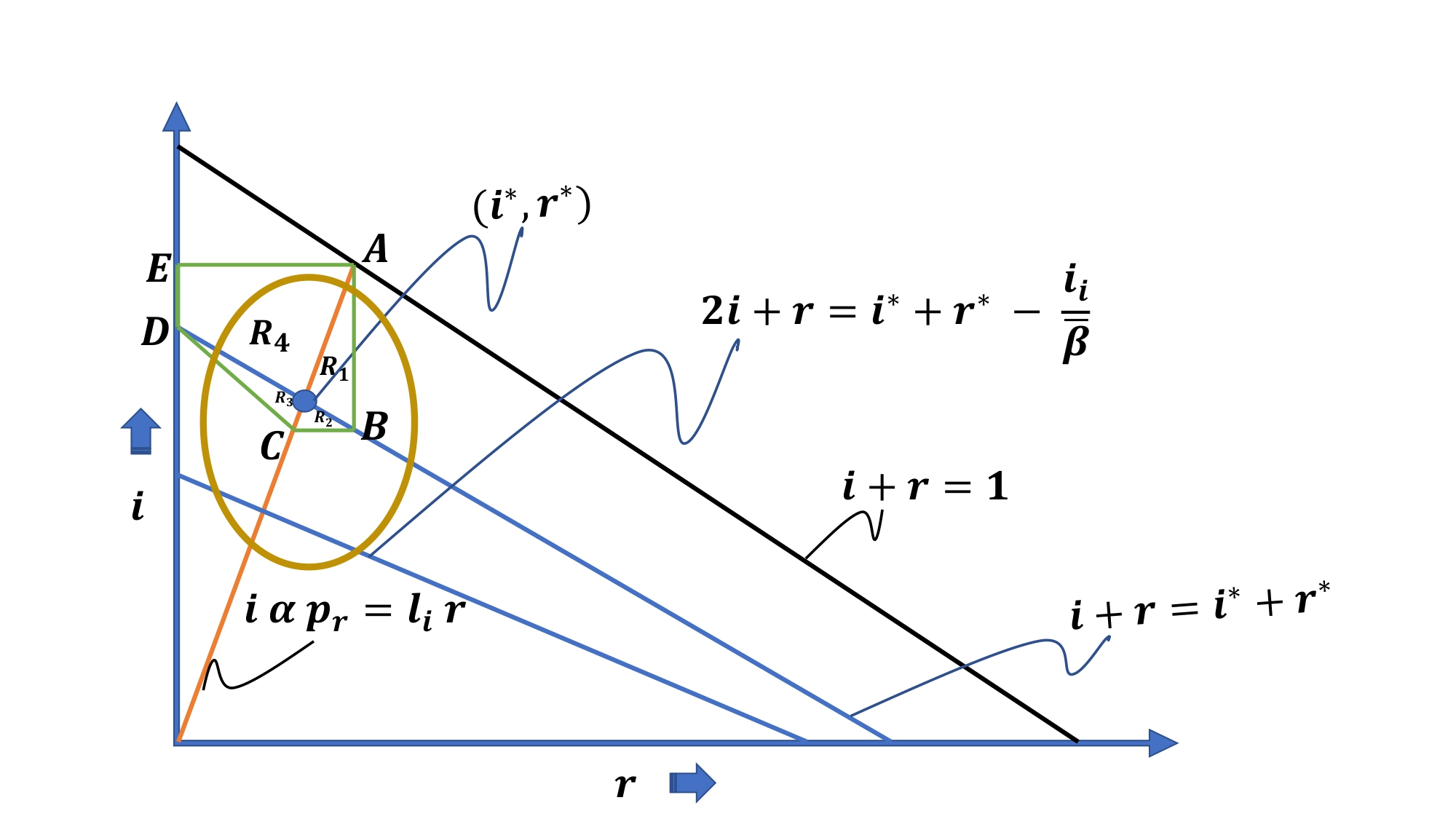}
  \vspace{-3mm}  \caption{Limit cycle when $\alpha + l_i < \bar \beta$,  $\bar \beta > 1 + \frac{l_i}{\alpha p_r}$}
   \vspace{-1mm}
\end{figure}}
That means $B^{(2)} /{\cal R}  $ does not contain any limit cycle, because any such limit cycle should intersect\footnote{This is not possible without touching  ${\cal R}$ (as in Fig. \ref{fig1}, the region ${\cal R}$ touches the boundary line $\{i+r=1\}$ at A).} the region ${\cal R}$ as it has to enclose $(i^*,r^*)$, which contradicts the fact that all points on or inside ${\cal R}$  are attracted toward $(i^*, r^*)$.  This completes the proof.
\eop

\noindent $\textbf{Proof of Theorem \ref{thm_Increased_interest}:}$  For ODE 
\eqref{Eqn_linear_influence_ode}, again $(0,0)$ is an equilibrium point and the other equilibrium points\footnote{    The numerator of $\beta(i) (1-i-r) -\alpha$, with  $r = \alpha p_r /l_i i$, is the quadratic function, $c_ai^2 + c_b i + c_c$. } are: 

\vspace{-3mm}
{\small \begin{eqnarray}
\label{Eqn_thm3_eqpt}
i^*= \frac{-c_b \pm  \sqrt{(c_b)^2 - 4c_a c_c}}{2c_a},
\ r^* = \frac{\alpha p_r}{l_i} i^*.
\end{eqnarray}}
We now identify the ones in $\B$ (equivalently $i^* \in [0, 1/c_i^r]$, for {\small$c_i^r := 1+ \alpha p_r/l_i$}) and further identify if they are attractors. Consider the corresponding Jacobian matrices:

\vspace{-3mm}
{\small$$
J (i^*, r^*) = \begin{bmatrix}
  i^*(2 i^* c_a + c_b) & -i^* \bar\beta   \\
\alpha p_r & -l_i
\end{bmatrix}
\mbox{, }
J (0, 0) = \begin{bmatrix}
   \beta_0 - {\alpha} & 0 \\
\alpha p_r & -l_i
\end{bmatrix}.$$}

\TR{}{

{\bf Case I, when {\small$\beta_0 - {\alpha}<0$}:} 
Begin with sub-case ${\bar w} <0$ and define, $L(i) = \eta (pi+q)$.  If   $
 {\bar \beta} + {\bar w} \frac{ L(i) }{a L(i) + 1} < 0$, then 
  from ODE \eqref{Eqn_linear_influence_ode}, we 
 clearly have:
 $$
 \frac{di} {dt}  < - i \alpha
 < -i \min\{\alpha,  \alpha - \beta_0\}. $$
 On the other hand, when $
 {\bar \beta} + {\bar w} \frac{ L(i) }{a L(i) + 1} \ge  0$, we have
 \begin{eqnarray*}
 \frac{di} {dt} 
 &\le& i  \left (   {\bar \beta} + {\bar w} \frac{ L(i) }{a L(i) + 1}   - \alpha \right ) \\
 &\stackrel{a}{\le}& i (\beta_0 - \alpha)  \ \
 <  \  -i \min\{\alpha,  \alpha - \beta_0\}. 
 \end{eqnarray*}and the inequality  `$a$' holds because  the
 function 
 $i \mapsto L(i)/ (a L(i) + 1)$ 
    is increasing and $\bar{w} < 0$.
    Thus in all, with ${\bar w} < 0$, we have
  \begin{eqnarray*}
 \frac{di} {dt} 
 &\le &   -i \min\{\alpha,  \alpha - \beta_0\}. \end{eqnarray*}
 By   standard results in ODEs, we have  that (also given as \cite[Lemma 6]{agarwal2021new}),  $i(t) \le i_0 \exp\left ({- \min\{\alpha,  \alpha - \beta_0\} t} \right ) \to 0 $. 
 Also $i(t) \ge 0$ by Lemma \ref{lemma1}. Hence, $i(t) \to 0$ and  it is simple to show  that $(0,0)$ is a global attractor on $\B$, when ${\bar w} < 0$ (see r-ODE in \eqref{Eqn_combined_ODE}).

 Now, consider the sub-case with  $\beta_0 < \alpha$,
and ${\bar w} > 0$, but $\beta a + {\bar w} - \alpha a < 0$; then 
   one can upper bound the $i$-ODE as: 
$$
\frac{di}{dt} \le 
i \left ( \beta + \frac{{\bar w} }{ a }- \alpha \right )
= \frac{i}{a} (\beta a + {\bar w} - \alpha a ) < 0,
$$ and hence $(i(t), r(t)) \to (0, 0) $, as before,  for any i.c.

\ignore{
The leftover case, when $\beta_0 < \alpha$, $\bar w > 0$ and ${\bar w} > 0$, but $\beta a + {\bar w} - \alpha a > 0$, needs to be completed with probably the following:
1) when $c_b< 0$, then nothing to prove (as all the roots are negative), 2) $c_b >0$, then we can show that $c_b^2 < 4c_a c_c$.
}}

\TR{By using similar logic as before, $(0,0)$ is a global attractor in this scenario (details in \cite{TR}). Now consider the case    with $\beta_0 - {\alpha}>0$.}{

{\bf Case II,  with $\beta_0 - {\alpha}>0$:}} Define the function $$p_1 (i) = c_a i^2 + c_b i + c_c, \mbox{ for } i \in [0,1/c_i^r].$$ Then, $p_1 (0) = (\beta_0 - {\alpha})(a \bar{\eta}q +1) > 0$ and, $p_1(1/c_i^r) = -\alpha < 0$. Thus, by continuity, there exists a unique zero of $p_1(i) \in (0, 1/c_i^r]$, and is one of the two values in \eqref{Eqn_thm3_eqpt}.

\TR{One can show that the root inside $[0, 1/c_i^r]$, is the one with the negative sign using simple algebra (details in \cite{TR}).}{
  The square root term will be $\sqrt{c_b^2 - 4 c_a c_c} > \mbox{ or } < |c_b|$, when $c_a$ is positive or negative, respectively. In both   cases, the negative sign of square root term   leads to   the required root, irrespective of $c_b > 0$ or $c_b < 0$. Thus,
  
\vspace{-4mm}
    {\small 
    \begin{eqnarray*}
  i^*= \frac{-c_b -  \sqrt{(c_b)^2 - 4c_a c_c}}{2c_a},
\ r^* = \frac{\alpha p_r}{l_i} i^*.
\end{eqnarray*}}

}
The corresponding Jacobian matrix $J(i^*, r^*)$ is negative definite, making $(i^*, r^*)$ an LA  because

\vspace{-4mm}
{\small $$
2i^* c_a + c_b = - \left (c_b + \sqrt{c_b^2 - 4 c_a c_c} \right ) + c_b = -\sqrt{c_b^2 - 4 c_a c_c} < 0.
$$}The proof now follows by Lemma \ref{lemma1}  and PB theorem.
\eop

\noindent  \textbf{ Proof of Theorem \ref{sir_linear_beta}:} 
For   ODE \eqref{Eqn_linear_w}  using Jacobian matrix, $(0,0)$ is an LA iff $ \bar\beta < \alpha $. The possible nonzero equilibrium points are  given by:
 
  \vspace{-4mm}
    {\small 
    \begin{eqnarray}
    \label{Eqn_i*_r*}
i^* = \frac{(u -\bar\beta c^r_i) \pm  \sqrt{(\bar\beta c^r_i - u )^2 + 4u c^r_i(\bar\beta - \alpha) }}{2 u c^r_i},  \ r^* = \frac{ i^* \alpha p_r}{l_i}.
\end{eqnarray}}
When ${\bar \beta}  < \alpha$ and $u < \alpha c_i^r$, using simple algebra\TR{}{\footnote{ When $u < {\bar \beta} c_i^r$, both the roots are negative. When $  {\bar \beta} c_i^r < u <  \alpha  c_i^r $, 
then the discriminant term $= (\bar \beta c_i^r +u)^2 - 4 u \alpha c_i^r < 0$, so roots are not real.
}} one can show that  these non-zero equilibrium points will not  lie in interval $[0,1/c^r_i]$. \TR{Thus part (i) follows as before (see \cite{TR}).} {
Thus, by PB theorem the limit set with respect to any i.c.,  inside $B^{2}$, is $\{(0,0)\}$. Further since $\{(0,0)\}$ is an LA, the limit point becomes the limit and hence the global attractor.}

 When $\bar{\beta} > \alpha$, 
\TR{after some algebraic manipulations (details are in \cite{TR}), one can show that only the root with $+$ve sign in  \eqref{Eqn_i*_r*}}{
Define the function $$p_2 (i) =  \left ( \bar\beta +  u i \right ) (1-c_i^r i) -  \alpha, \mbox{ for } i \in [0,1/c_i^r].$$
 $ p_2(0) = (\bar{\beta} - {\alpha}) >0$ and  $ p_2(1/c_i^r) =  - {\alpha} <0$. Thus by continuity, only one non-zero root given by \eqref{Eqn_i*_r*} lies in   $[0, 1/c_i^r]$.

On following similar arguments as before,  $\sqrt{(\bar\beta c^r_i - u )^2 + 4u c^r_i(\bar\beta - \alpha) } > \mbox{ or } < |u-\bar{\beta}c^r_i |$, if $u$ is positive or negative respectively. So the only root with positive sign of square root will work, irrespective of the sign of $(u-\bar{\beta}c^r_i)$.

{\ignore{
{\bf Sub-case $u>0$:} In this case, we have
$2u c_i^r > 0 \mbox{ and } \sqrt{(\bar\beta c^r_i - u )^2 + 4u c^r_i(\bar\beta - \alpha) } > |u-\bar{\beta}c^r_i |$. So the only root with positive sign of square root will work, irrespective of the sign of $(u-\bar{\beta}c^r_i)$.

{\bf Sub-case $u<0$:} In this case, we have
$2u c_i^r < 0 \mbox{ and } \sqrt{(\bar\beta c^r_i - u )^2 + 4u c^r_i(\bar\beta - \alpha) } < |u-\bar{\beta}c^r_i |$. So again, the only root with positive sign of square root will work.}
}
Hence, in both cases, 

\vspace{-4mm}
    {\small 
    \begin{eqnarray*}
   i^* = \frac{(u -\bar\beta c^r_i) +  \sqrt{(\bar\beta c^r_i - u )^2 + 4u c^r_i(\bar\beta - \alpha) }}{2 u c^r_i},
   \   r^* = \frac{ i^* \alpha p_r}{l_i}. 
\end{eqnarray*}}

} 
   lies in $\B$. 
\ignore{
Also, on assuming \eqref{xy},  $(x,y) \to (0,0)$ and \\
\begin{eqnarray*}
    \frac{dx}{dt} \ = \  -i(\bar\beta(x+y) + u(x(i+i^* -1) +rx +i^* y)), \\   
    \  \  \   \ 
\frac{d y}{d t} \ = \  x\alpha p_r - y l_i.  
\end{eqnarray*}}
Now, at $ (i^*, r^*)$, the Jacobian matrix  is:

\vspace{-4mm}
{\small \begin{eqnarray}
J (i^*, r^*) = \begin{bmatrix}
    i^* (u(1-i^* - r^*)-\bar\beta - u i^*) & -i^* (\bar\beta + u  i^*)  \\
\alpha p_r & -l_i
\end{bmatrix}
\label{thm4_jacobin}.
\end{eqnarray}}

The product of eigenvalues of \eqref{thm4_jacobin} is greater than zero by Lemma \TR{3 in \cite{TR}}{\ref{lemma3}}, 
 and hence $(i^*,r^*)$ is either a source or sink. Thus it is a negative definite iff
$
u  (1-c_i^r i^*)  \le 
{\bar \beta} + u i^* + l_i /i^*,
$ otherwise, $(i^*, r^*)$ will be a source.
Thus part (ii).a follows by  PB theorem, as $(i^*, r^*)$ is now an LA. 

 For part (ii).b, $(i^*, r^*)$ is unstable. 
 %
 There cannot be a closed orbit (or limit cycle) that just encloses $(0,0)$ by the critical point criterion (e.g.,  \cite{palis2012geometric}), as $(0,0)$  is a saddle point. 
Further,   the  trajectory starting from $(r,0)$ in $\B$, will converge to $(0,0)$. We need to understand the remaining trajectories. The rest of the theorem follows from PB theorem, as $(i^*, r^*)$ is the only other equilibrium point, and it is not a saddle point; thus a trajectory in $\B$ approaches a  limit cycle  enclosing $(i^*, r^*)$, or approaches one of the two equilibrium  points.  \eop

{\ignore{
\newcommand{\LongProof}[1]{{\color{red}#1}}
\LongProof{
\section*{Proof of No Limit Cycle}
Step 0 - The trajectory is trapped inside $\{ (i, r) \in [0,1]^2: i+r \le 1 \}$
as in Theorem 4. 

Step1 - To construct a rectangle $R_1$ around $(i^*, r^*)$ which includes a part of line $2i+r = i^*+r^* - l_i / \beta$ and observe the vector field in pointing interior or on the boundary.  (Check a reference to show it gets trapped in such rectangles).

Step 2- By using Bendixson's criterion \cite[Theorem 3.3]{li1993bendixson} no limit cycle exists inside this smaller region.  Hence by Poincare Bendixson criterion \cite[Theorem 1.8]{palis2012geometric} trajectories from any point inside this region  converge towards $(i^*, r^*)$

Step-3 There can be no other region  inside $\{ (i, r) \in [0,1]^2: i+r \le 1 \}$ which contains limit cycle because, then a part of the LC should pass through region $R_1$ (as any LC should enclose a critical point by \cite[]{.}) which contradicts the fact that all points inside $R_1$ are attracted towards $(i^*, r^*)$. 

Step-4  Thus no limit cycle inside $R_0\{ (i, r) \in [0,1]^2: i+r \le 1 \}$ and the trajectory confined by \cite[Theorem ]{PB} all the points inside $R_0$ converge towards $(i^*, r^*)$. }
}}

{}

\TR{
 \begin{figure}[htbp] 
 \vspace{-3mm}
           \centering
           \begin{minipage}{4.2cm}
           \hspace{1mm}
           \includegraphics[width=0.9\textwidth]{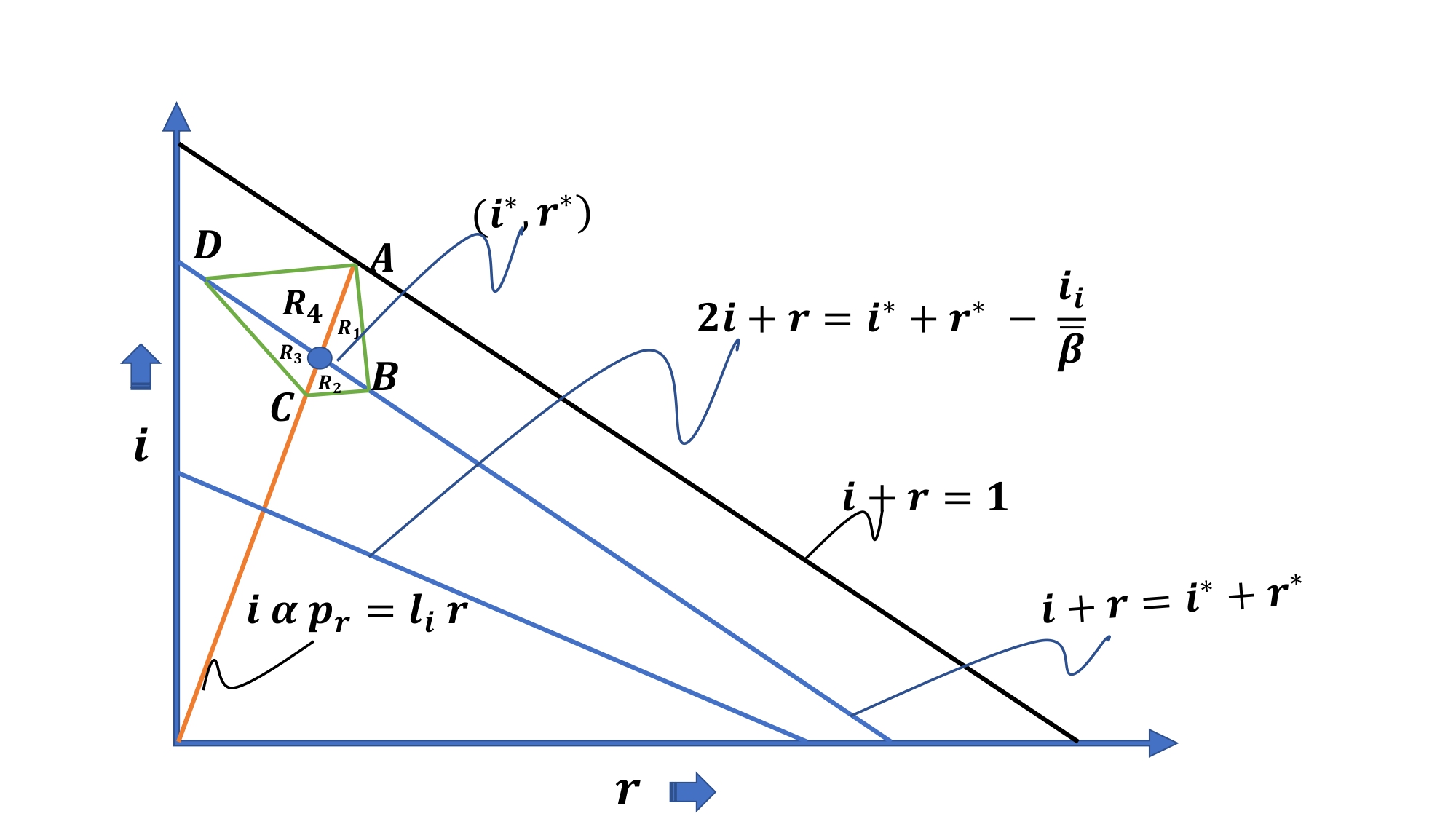}
           \end{minipage}
            \begin{minipage}{4.2cm}
            \hspace{1mm}
           \includegraphics[width=0.9\textwidth]{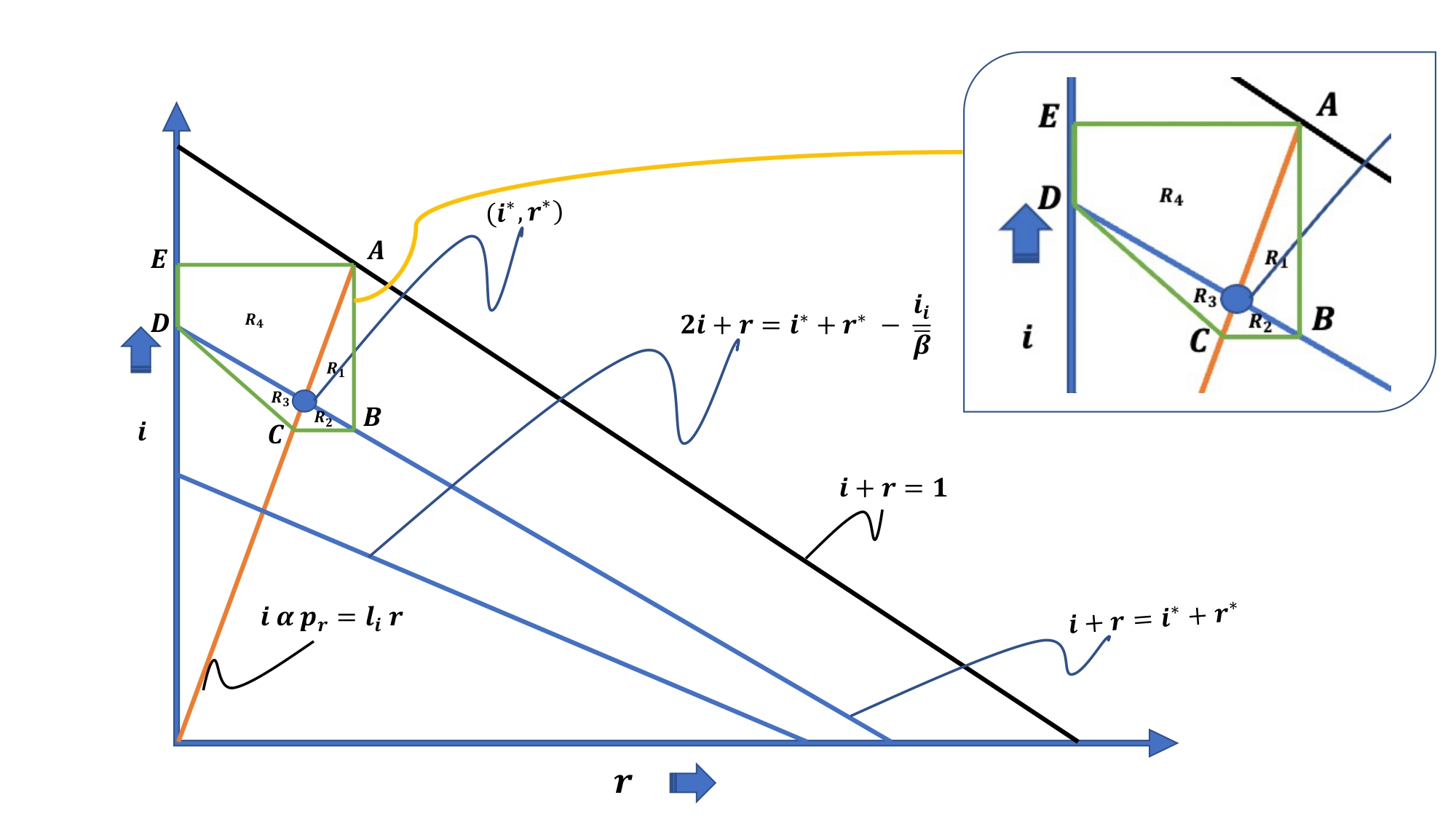}
             \end{minipage}  
    \caption{Left:   $\alpha + l_i < \bar \beta$,  $\bar \beta < 1 + \frac{l_i}{\alpha p_r}$, \ \ Right:   $\alpha + l_i < \bar \beta$, $\bar \beta > 1 + \frac{l_i}{\alpha p_r}$ \label{fig1}}
    \end{figure}} 
    {
  \begin{figure*}[htbp] 
 \vspace{-4mm}
           \centering
           \begin{minipage}{8.2cm}
           \hspace{-6mm}
           \includegraphics[width=1.1\textwidth]{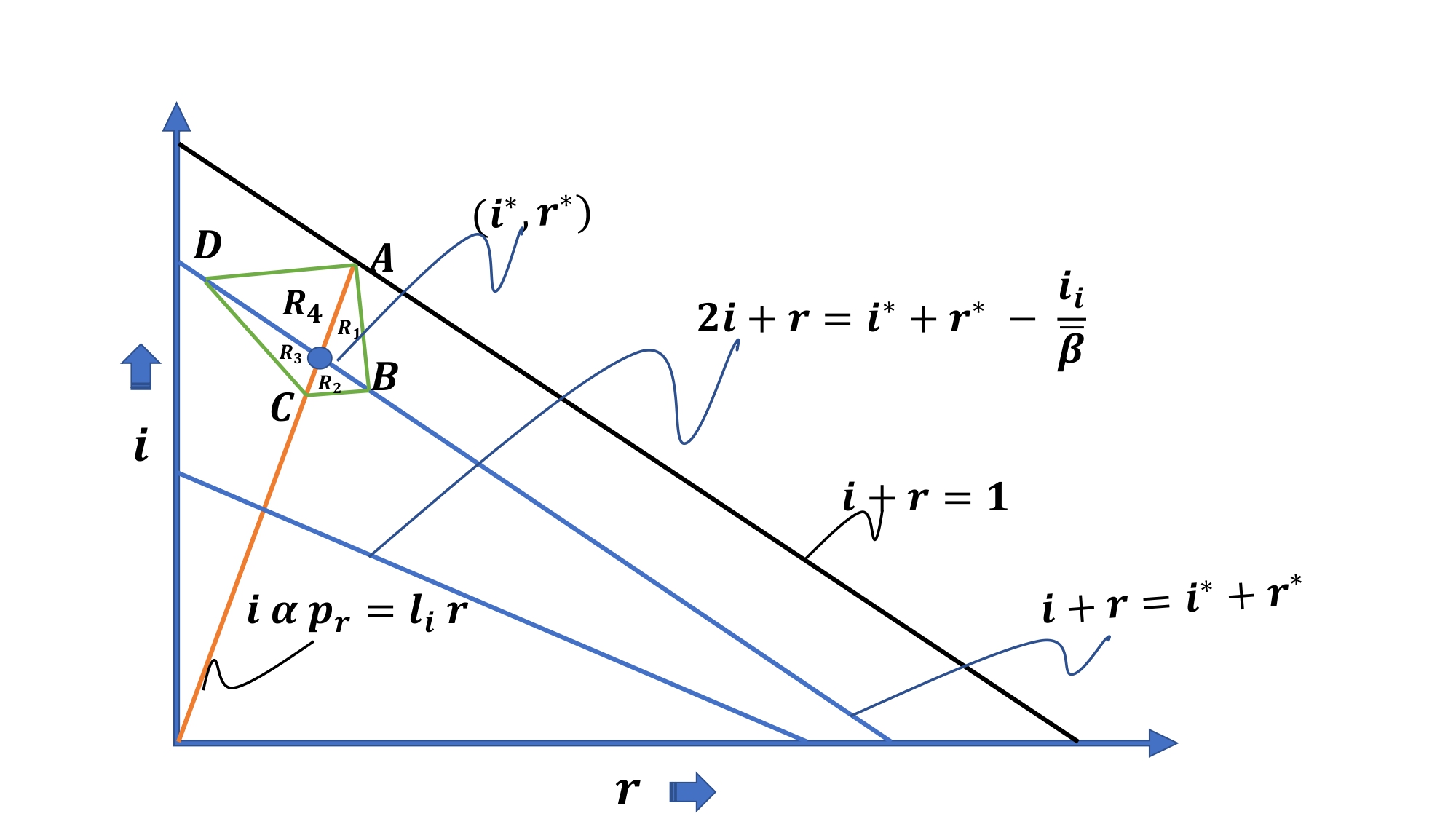}
           \end{minipage}
            \begin{minipage}{8.2cm}
            \hspace{-6mm}
           \includegraphics[width=1.1\textwidth]{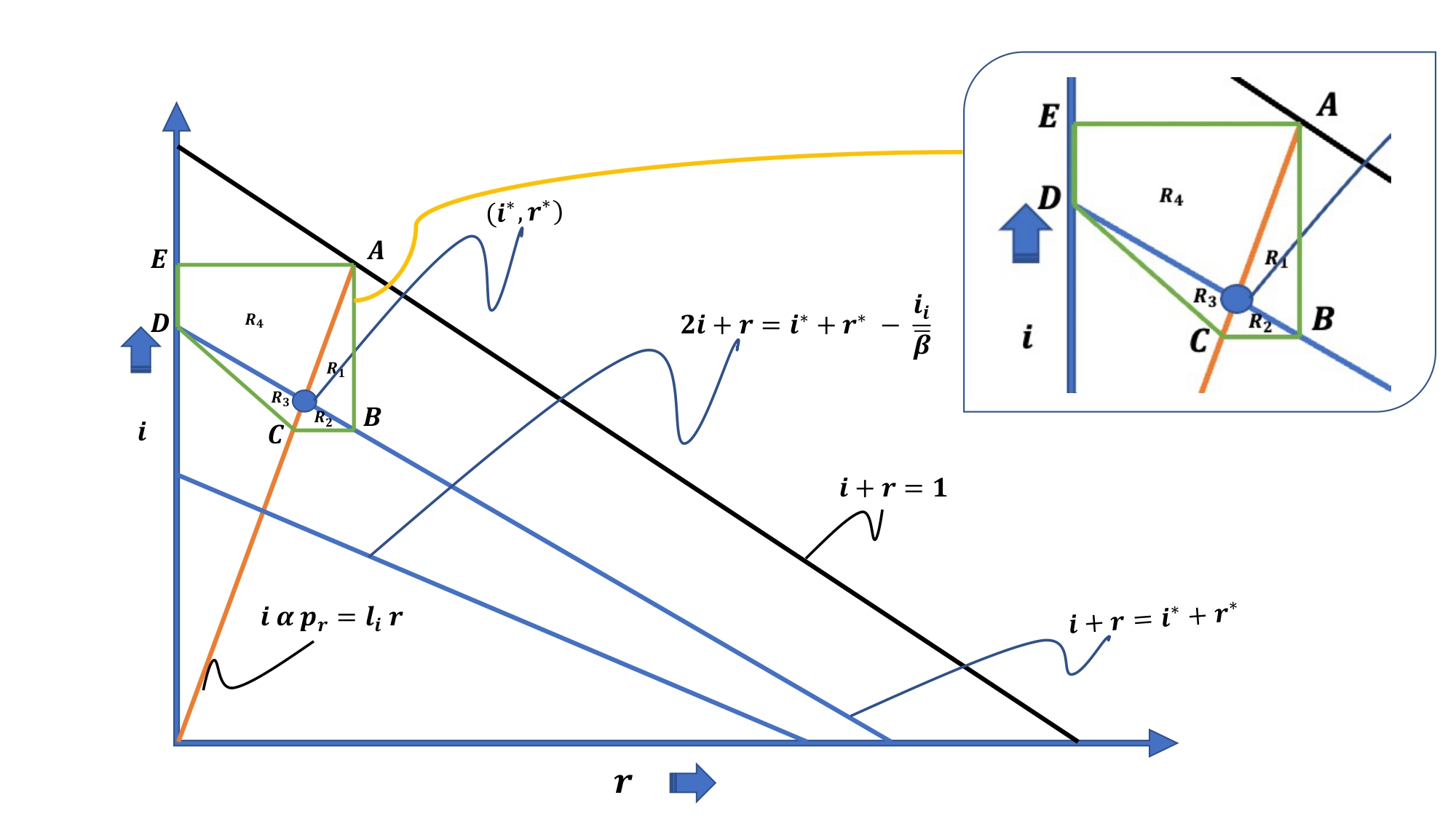}
             \end{minipage}  
    \caption{Left:   $\alpha + l_i < \bar \beta$,  $\bar \beta < 1 + \frac{l_i}{\alpha p_r}$, \ \ Right:   $\alpha + l_i < \bar \beta$, $\bar \beta > 1 + \frac{l_i}{\alpha p_r}$ \label{fig1}}  
        \end{figure*}
    }
\TR{}{
 \begin{lem}
 \label{lemma3}
 The product term of Jacobian matrix \eqref{thm4_jacobin} is positive.
 
 {\bf Proof:} The product term is
 \vspace{-1mm}
 {\small
 \begin{eqnarray*}
 i^*\left( -\left(\frac{u \alpha}{K} - K \right)l_i + K \alpha p_r \right) = \frac{i^* l_i}{K} \left( c_i^r K^2 - {u \alpha} \right),
 \end{eqnarray*}
 }
 
 where $ K := {\bar{\beta} + u i^*}$.
 Now, as $K = \frac{\alpha}{1-i^* - r^*} > 0$, we check the sign of $\left( c_i^r K^2 - {u \alpha} \right)$ term only. 
 So, consider
 \vspace{-1mm}
 {\small
 \begin{eqnarray*}
4 c_i^r ( c_i^r K^2 - {u \alpha}) \\
= \left(2 \bar{\beta}c_i^r + (u - \bar{\beta}c_i^r) + \sqrt{(u + \bar{\beta}c_i^r)^2 - 4\alpha u c_i^r } \right)^2  - 4\alpha u c_i^r \\ = 2(u+\bar{\beta}c_i^r)^2 - 8\alpha u c_i^r + 2(u+\bar{\beta}c_i^r)\sqrt{(u + \bar{\beta}c_i^r)^2 - 4\alpha u c_i^r },  \end{eqnarray*}}
is positive as $(u+\bar{\beta}c_i^r)^2 > 4\alpha u c_i^r$, because the square root term $(u+\bar{\beta}c_i^r)^2 - 4\alpha u c_i^r >0$ used for $(i^*, r^*)$ definition. Now as $c_i^r$ is positive, the product term is positive. This completes the proof. \eop
 \end{lem}  } 
 
 \begin{lem}
 \label{lemma2}
 Under the assumptions of Theorem \ref{thm_sir_attractor}.(ii), one can construct a closed, bounded region ${\cal R} \subset \B$, touching boundary $\{i+r=1\}\subset \partial (\B)$  at $A=(b_*,  1-b_* ) \in \partial {\cal R}$, where $b_*=l_i/(l_i+\alpha p_r)$, such that 

 \vspace{-4mm}
 {\small\begin{eqnarray}
 \label{Eqn_fg}
 \pdv{f}{i} + \pdv{g}{r} = -(\Bar{\beta}(2i+r-i^* - r^*) + l_i) \ne 0, 
 \mbox{ where } \nonumber\\
 f(i,r) = i \left ( \Bar{\beta}  (1-i-r) - \alpha \right ),   \mbox{ and, } g(i,r) = i \alpha p_r  - l_i  r,
 \end{eqnarray}}and such that the field  $(f,g)$ representing ODE \eqref{epi_ode} at the boundary   $\partial{\cal R} $ is  pointing inwards or onto the boundary    $\partial{\cal R} $.
 \end{lem}
 
 \TR{
 {\bf Proof:} 
The region ${\cal R}$ is constructed   in Figure \ref{fig1}, the rest of the details are in \cite{TR}. \eop}{}

\TR{}{\noindent {\bf Proof of Lemma \ref{lemma2}:} Here, we will construct a bounded region,  such that  the  field of the ODE on its boundary is either  inwards or on the boundary of that bounded region. Then, similar arguments will follow as lemma \ref{lemma1}. 
Further, the region can be constructed in such a way that the line given by 
\begin{equation}
\Bar{\beta}(2i+r-i^* - r^*) + l_i = 0,
 \label{Eqn_line}    
\end{equation}
is outside the region ${\cal R}$ (see Fig. \ref{fig1}). 
We will have three cases the following:

 \subsection{\bf{Consider the case when $\alpha + l_i < \bar \beta$ and $\bar \beta < 1 + \frac{l_i}{\alpha p_r}$ }} At the corner point  $A = \left( \frac{\alpha p_r}{\alpha p_r + l_i}, \frac{l_i}{\alpha p_r + l_i} \right )$, the flow is downwards (along line $AB$) as $dr/dt = 0$  and $di/dt < 0 $. Similarly, at corner point  $B$, flow is towards left (along line $BC$) as on line $i+r = i^* +r^*$,   ${di}/{dt} = 0 $  and $dr/dt < 0 $ is decreasing (as  $i \alpha p_r < r l_i$). Similar arguments will follow for corner points $C$ and $D$.  Now, consider any point  $(r,i) \mbox{ on } AB$, flow will be inwards as $di/dt <0$ and $dr/dt <0$ are decreasing  (because $(r,i) \in \{ (r,i) ; i+r > i^* + r^* \mbox{ and } i \alpha p_r < r l_i\}$). Similarly, for any point $(r,i)$ on boundary $BC$, have flow upwards as $i$ is increasing below line $(i+r = i^* +r^*)$ and  $r$ is decreasing below line $(i \alpha p_r = r l_i)$. Similar arguments will follow for boundary lines $CD$ and $DA$. 

 Furthermore, observe that the above arguments only require that $AB$ is vertical ($A$ touching the lines $i+r=1$ and  $\alpha p_r i = l_i r$),  both $AD$ and $BC$ are horizontal and CD is then joining the points. Thus one can shift the line $BC$ upward if required to ensure, $C$ does not touch the line given in \eqref{Eqn_line} (observe that $2i+l_i/\beta = i$ only when ${\bar \beta}  <0$, which is not true here, thus there is a gap between lines $i+r = i^*+r^* $ and the line given in \eqref{Eqn_line}). 
 
  \subsection{\bf{Consider the case when $\alpha + l_i < \bar \beta$ and $\bar \beta > 1 + \frac{l_i}{\alpha p_r}$ }}
    
  On following the similar arguments as above, one can prove the trajectory will be confined in the $ABCDE$ region when started from anywhere inside the same region. Further again, the lines $2i+l_i/\beta + r = i^*+r^*$   and $i+r = i^*+r^*$ have a gap in between because $\beta  > 0.$

  \subsection{\bf{Consider the case when $\alpha + l_i > \bar \beta$ }}

 In this scenario, by the Bendixson criterion, we don't have any limit cycle in the $\B$ region.  
 And further a bounded region as in previous cases can be constructed, 
 
 This completes the proof by applying Nagumo’s theorem \cite{ghorbal2022characterizing}. 
 \eop
}

}

\end{document}